\documentclass[a4paper,10pt]{amsart}

        \title{Commuting homotopy limits and smash products}
       \author{Wolfgang L\"uck}
       \author{Holger Reich}
       \author{Marco Varisco}
         \date{August 21, 2003}
      \address{Westf\"alische Wilhelms-Universit\"at M\"unster, 
               Mathematisches Institut. 
               Einsteinstr.~62, 
               D-48149 M\"unster, Germany}
        \email{lueck@math.uni-muenster.de}
      \urladdr{http://www.math.uni-muenster.de/u/lueck}
        \email{reichh@math.uni-muenster.de}
      \urladdr{http://www.math.uni-muenster.de/u/reichh}
        \email{varisco@math.uni-muenster.de}
      \urladdr{http://www.math.uni-muenster.de/u/varisco}

\usepackage[arrow,curve,matrix,tips]{xy} \SelectTips{eu}{10} \UseTips

\newcommand{\calrf}{\mathcal{RF}}

\newcommand{\calf}{\mathcal{F}}
\newcommand{\calg}{\mathcal{G}}  

\newcommand{\cali}{\mathcal{I}}

\newcommand{\caln}{\mathcal{N}}

\newcommand{\calr}{\mathcal{R}}

\newcommand{\bbN}{\mathbb{N}}

\newcommand{\bbQ}{\mathbb{Q}} 
\newcommand{\bbR}{\mathbb{R}}

\newcommand{\bbZ}{\mathbb{Z}} 

\newcommand{\bfA}{\mathbf{A}}    \newcommand{\bfa}{\mathbf{a}}
    \newcommand{\bfb}{\mathbf{b}}

\newcommand{\bfE}{\mathbf{E}}    
\newcommand{\bfF}{\mathbf{F}}    \newcommand{\bff}{\mathbf{f}}
    \newcommand{\bfg}{\mathbf{g}}

\newcommand{\bfR}{\mathbf{R}}    \newcommand{\bfr}{\mathbf{r}}
    
    \newcommand{\bft}{\mathbf{t}}

\DeclareMathAlphabet{\matheurm}{U}{eur}{m}{n}

\newcommand{\SPACES}{\matheurm{SPACES}}
\newcommand{\SPECTRA}{\matheurm{SPECTRA}}
\newcommand{\MODULES}{\matheurm{MODULES}}
\newcommand{\GROUPOIDS}{\matheurm{GROUPOIDS}}
\newcommand{\Or}{\matheurm{Or}}
\newcommand{\Sub}{\matheurm{Sub}}

\newcommand{\sma}{\wedge}
\newcommand{\ot}{\leftarrow}
\newcommand{\EGF}[2]{E{#1}(#2)}

\newcommand{\Ann}{\operatorname{Ann}}
\newcommand{\aut}{\operatorname{aut}}
\newcommand{\ch}{\operatorname{ch}}
\newcommand{\colim}{\operatorname{colim}}
\newcommand{\conhom}{\operatorname{conhom}}
\newcommand{\Ext}{\operatorname{Ext}}
\newcommand{\hocolim}{\operatorname{hocolim}}
\newcommand{\holim}{\operatorname{holim}}
\newcommand{\hopb}[3]{\operatorname{hopb}(#1\!\to\!#2\!\ot\!#3)}
\newcommand{\hopo}[3]{\operatorname{hopo}(#1\!\ot\!#2\!\to\!#3)}
\newcommand{\id}{\operatorname{id}}
\newcommand{\im}{\operatorname{im}}
\newcommand{\Inn}{\operatorname{Inn}}
\newcommand{\map}{\operatorname{map}}
\newcommand{\mor}{\operatorname{mor}}
\newcommand{\obj}{\operatorname{obj}}
\newcommand{\op}{\operatorname{op}}
\newcommand{\pr}{\operatorname{pr}}
\newcommand{\pt}{\operatorname{pt}}
\newcommand{\tors}{\operatorname{tors}}
\newcommand{\Tor}{\operatorname{Tor}}

\newcommand{\comsquare}[8]{%
\[
\xymatrix{
  #1 \ar[d]^-{#4} \ar[r]^-{#2} & #3 \ar[d]^-{#5} \\
  #6 \ar[r]^-{#7} & #8
}
\]}

\theoremstyle{plain}
  \newtheorem{theorem}               {Theorem}     [section]
  \newtheorem{lemma}       [theorem] {Lemma}
  \newtheorem{corollary}   [theorem] {Corollary}
  \newtheorem{proposition} [theorem] {Proposition}
  \newtheorem{addendum}    [theorem] {Addendum}
  \newtheorem{question}    [theorem] {Question}

\theoremstyle{definition}
  \newtheorem{definition}  [theorem] {Definition}
  \newtheorem{example}     [theorem] {Example}

\theoremstyle{remark}
  \newtheorem{remark}      [theorem] {Remark}

\makeatletter\let\c@equation=\c@theorem\makeatother

\begin{document}

\typeout{--------------------  Abstract  --------------------}

\begin{abstract}
  In general the processes of taking a homotopy inverse limit of a
  diagram of spectra and smashing spectra with a fixed space do not
  commute. In this paper we investigate under what additional
  assumptions these two processes do commute. In fact we deal with an
  equivariant generalization that involves spectra and smash products
  over the orbit category of a discrete group.  Such a situation
  naturally occurs if one studies the equivariant homology theory
  associated to topological cyclic homology.  The main theorem of this
  paper will play a role in the generalization of the results obtained
  by B\"okstedt, Hsiang and Madsen about the algebraic $K$-theory
  Novikov Conjecture to the assembly map for the family of virtually
  cyclic subgroups.
\end{abstract}

\maketitle

\begin{flushright}
  \textit{to Hyman Bass on his seventieth birthday}
\end{flushright}
\medskip

\typeout{--------------------  Introduction  --------------------}

\section{Introduction}
\label{sec: Introduction}

Let $\cali$ be a small category and let $\bfE: \cali \to \Omega \mbox{-} \SPECTRA$
be a contravariant functor, then for any space~$X$ there is a natural map
\[
X_+ \sma \holim_{\cali} \bfE \to \holim_{\cali} X_+ \sma \bfE.
\]
In general this map is not a weak equivalence. 
A special case of our main result says that the map is an equivalence
if we assume that:
\begin{enumerate}
\item
The spectra $\bfE (c)$ for objects $c$ in $\cali$ are uniformly bounded below.
\item
The homology of $X$ is degreewise finitely generated as an abelian group.
\item
The category $\cali$ admits a finite dimensional contravariant
classifying space (compare Definition~\ref{def: Ecalc}).
\end{enumerate}

In fact we will generalize the question in several directions.  On the
one hand we will work throughout in an equivariant setting. On the
other hand we will investigate general mapping space constructions, of
which homotopy limits are just special cases.  Moreover we will also
try to obtain weaker assumptions that still suffice to conclude that
the map induces an isomorphism on rationalized homotopy groups. And
finally we will discuss a chain complex analogue of our result.

In order to state our main result we introduce some more notation.
Let $G$ be a discrete group.
Let $\Or(G)$ be the \emph{orbit category},
whose objects are the homogeneous left $G$-spaces $G/H$ and
whose morphisms are $G$-maps.
Every left $G$-space~$X$, in particular every $G$-CW-complex
(compare Section~\ref{sec: The homology of classifying spaces for families}),
gives rise to 
a contravariant  $\Or(G)$-space, i.e.~a contravariant  functor from the orbit category to spaces,
by sending $G/H$ to $X^H = \map_G ( G/H , X)$.
If we furthermore have a covariant functor $\bfE$ from $\Or(G)$ to $\SPECTRA$ (from now on called 
an $\Or (G)$-spectrum) then we can form the smash 
product over the orbit category
\[
X_+ \sma_{\Or (G)} \bfE 
\]
to obtain an ordinary spectrum 
(compare Section~\ref{sec: Basics about modules and spectra over a category}). 
This is an important construction since 
sending $X$ to  $\pi_{\ast} ( X_+ \sma_{\Or (G)} \bfE )$ yields an equivariant homology theory.

Our basic conventions concerning the categories $\SPECTRA$ and $\Omega$-$\SPECTRA$ can be found
in Section~\ref{sec: Basics about spectra}. In particular we will denote by 
\[
\bfR: \SPECTRA \to \Omega \mbox{-} \SPECTRA
\]
a \emph{fibrant
replacement} functor that replaces a spectrum by a weakly equivalent $\Omega$-spectrum, and $\bfr \colon \bfE \to \bfR ( \bfE )$
will denote the natural equivalence from a spectrum to its replacement,
compare~\eqref{bfr}.

The \emph{homotopy limit} of a contravariant functor $\bfE : \cali \to \SPECTRA$ is defined as 
\[
\holim_{\cali} \bfE = \map_{\cali} ( E \cali_+ , \bfR \circ \bfE ) ,
\]
see Definition~\ref{def:holim}.
(If $\bfE$ takes already values in $\Omega\text{-}\SPECTRA$, then one can work with $\bfE$ instead of 
$\bfR \circ \bfE$, but in general the replacement is crucial.)
The $\cali$-mapping spectrum construction is explained in Section~\ref{sec: Basics about modules and spectra over a category} and 
$E \cali$ is the contravariant classifying space of the category $\cali$ (see Definition~\ref{def: Ecalc}).
In particular $E \cali$ is a 
contravariant $\cali$-$CW$-complex, i.e.~a contravariant functor from $\cali$ to spaces which is built up
out of free cells (see Definition~\ref{def:  cali-CW-complexes}). 
Even if one is eventually only interested in 
homotopy limits it is important for the proof of our main result to deal with other $\cali$-$CW$-complexes in
place of $E \cali$, because  the proof will proceed via induction over the skeleta.

We will now state the question we want to investigate in full generality.
Let $\cali$ be a small category and $Y$ a contravariant $\cali$-$CW$-complex.
Furthermore let
\[
\bfE \colon \cali^{\op} \times \Or (G) \to \Omega \mbox{-}\SPECTRA
\]
be a functor and $X$ a $G$-$CW$-complex.
There are maps of spectra $\bft'$ and~$\bft''$, which are natural in~$X$, $Y$ and~$\bfE$. They are given as follows:
\[
\bft'\colon X_+ \sma_{\Or(G)} \map_{\cali}(Y_+,\bfE) 
\to
\map_{\cali}(Y_+,X_+ \sma_{\Or(G)} \bfE)
\]
sends $x \wedge \phi$ to the map from $Y$ to  $X_+ \wedge \bfE$ which is given by
$y \mapsto x \wedge \phi(y)$  and
\[
\bft''\colon \map_{\cali}(Y_+,X_+ \sma_{\Or(G)} \bfE) 
\to
\map_{\cali}\left(Y_+,\bfR\left(X_+ \sma_{\Or(G)} \bfE\right)\right)
\]
is induced by the map $\bfr\left(X_+ \sma_{\Or(G)}  \bfE\right)$, the weak equivalence to the fibrant replacement.
\begin{question} \label{question}
Under what assumptions is the map  
\[
\bft\colon  X_+ \sma_{\Or(G)} \map_{\cali}(Y_+,\bfE) 
\to 
\map_{\cali}\left(Y_+,\bfR\left(X_+ \sma_{\Or(G)} \bfE\right)\right),
\]
defined as the composition $\bft'' \circ \bft'$, a weak equivalence?
\end{question}

We still need to explain a mild condition that we will always impose on our functor $\bfE$.
Let $\GROUPOIDS$ denote the category of small groupoids and functors between them.
For a $G$-set $S$ let 
$\calg^G( S)$ denote the \emph{transport groupoid}, whose set of objects is $S$ and for which
the set of morphisms from $s_1$ to $s_2$ consists of $\{ g \in G | gs_1 = s_2 \}$.
Composition of morphisms is induced from the group structure on $G$. We obtain a 
covariant functor $\calg^G \colon \Or (G) \to \GROUPOIDS$ by assigning $\calg^G (G/H)$ to $G/H$.
We call a functor $\bfF \colon \GROUPOIDS \to \SPECTRA$ a \emph{homotopy functor} if it sends
an equivalence of groupoids to a weak equivalence of spectra. 
It seems that all interesting examples of $\Or (G)$-spectra can be obtained by composing 
a suitable homotopy functor $\bfF: \GROUPOIDS \to \SPECTRA$ with the functor $\calg^G$.
For a subgroup $H \subset G$ we denote by 
$Z_GH = \{g \in G\mid gh = hg \text{ for all } h \in H\}$
its \emph{centralizer}.

\begin{theorem}\label{the: commuting wedge and homotopy limit}
Consider the following data:
\begin{itemize}
\item
A functor $\bfE \colon \cali^{\op} \times \Or(G) \to \Omega\text{-}\SPECTRA$
which factorizes as
\[
\bfE = \bfF \circ (\id \times \calg^G),
\]
where $\bfF$ is such that
for each object $c$ in~$\cali$
the functor $\bfF(c, - )$ is a homotopy functor.
\item
A contravariant $\cali$-$CW$ complex~$Y$.
\item
A $G$-CW complex~$X$.
\end{itemize}
Suppose that there are numbers $d$, $n$ and~$N \in \bbZ$ with~$d \geq 0$
such that the following conditions are satisfied:
\begin{enumerate}
\item[(A)]
$Y$ is $d$-dimensional.
\item[(B)]
The spectra $\bfE(c , G/H)$ are uniformly $(N-1)$-connected,
i.e.~for all objects $c$ in~$\cali$ and all orbits~$G/H$
we have $\pi_q (\bfE(c,G/H)) = 0 $ for~$q < N$.
\item[(C)]
The $G$-CW complex~$X$ has only finite isotropy groups
and only finitely many orbit types.
\item[(D)]
For each finite subgroup $H \subset G$ and every $p \leq n + d - N$
the homology group
\[
H_p ( Z_G H \backslash X^H ; \bbZ )
\]
is a finitely generated $\bbZ$-module.
\end{enumerate}
Then for each $p \leq n$ the map
\[
\pi_p (\bft) \colon  \pi_p \left( X_+ \sma_{\Or(G)} \map_{\cali}(Y_+,\bfE) \right)
\to
\pi_p \left( \map_{\cali}\left(Y_+,\bfR\left(X_+ \sma_{\Or(G)} \bfE\right)\right)\right)
\]
is an isomorphism.
\end{theorem}

This theorem will be proven in Section~\ref{sec: Commuting homotopy limits and smash products}.
In  Example~\ref{exa: necessary finiteness conditions on X} 
we discuss that none of the assumptions can be dropped.
A chain complex version of Theorem~\ref{the: commuting wedge and homotopy limit} will be 
discussed in Section~\ref{sec: The linear version}. 

If one is only interested in the question whether the map $\bft$
induces an isomorphism on rational homotopy groups one can obtain
results where assumption (D) is replaced by a weaker assumption.  A
very naive guess here would be to require the rationalized homology groups to be
finitely generated as $\bbQ$-vector spaces.  But this is too weak, compare Remark~\ref{rem: Q-finitely generated too weak}.
In order to formulate the weakest possible assumption and strongest conclusion
we introduce the notions of almost finitely generated modules and
almost isomorphisms, compare Section~\ref{sec: Almost trivial and almost finitely generated modules}.

A $\bbZ$-module $M$ is \emph{almost trivial}
if there is an element $r \in \bbZ, r \not=0$ such that $rm = 0$ holds for
all $m \in M$. A $\bbZ$-module $M$ is \emph{almost finitely generated} if
$M/\tors(M)$ is a finitely generated $\bbZ$-module and $\tors(M)$ is almost trivial.
A $\bbZ$-homomorphism is an \emph{almost isomorphism} if its kernel and cokernel are
almost trivial. An almost isomorphism becomes an isomorphism after
rationalization.

\begin{addendum}  \label{addendum}
If we weaken assumption (D) in Theorem~\ref{the: commuting wedge and homotopy limit} 
by requiring the homology groups to be almost finitely generated
$\bbZ$-modules, then we obtain the conclusion that the map
$\pi_p ( \bft )$ is an almost isomorphism and in particular a rational isomorphism for
all $p \leq n$.
\end{addendum}

As already explained we obtain results about homotopy limits by
specializing to classifying spaces.  For the convenience of the reader
we formulate the corresponding result explicitly.
\begin{corollary}
If $\cali$ admits a finite dimensional model for its contravariant classifying space $E \cali$, then
there is a  natural equivalence
\[
X_+ \sma_{\Or (G)} \holim_{\cali} \bfE \to \holim_{\cali} X_+ \sma_{\Or (G)} \bfE
\]
for a functor $\bfE: \cali^{\op} \times \Omega \mbox{-} \SPECTRA$ as in Theorem~\ref{the: commuting wedge and homotopy limit}
and a $G$-CW complex $X$ provided that:
\begin{enumerate}
\item[(B)]
The spectra $\bfE ( c , G/H )$ are uniformly bounded below.
\item[(C)]
The $G$-CW complex $X$ has only finite isotropy groups and only finitely many orbit types.
\item[(D)]
For every finite subgroup $H \subset G$ and every $n$ the homology group
\[
H_n ( Z_G H \backslash X^H ; \bbZ )
\]
is a finitely generated $\bbZ$-module.
\end{enumerate}
\end{corollary}

Of course there is again a version with an ``almost finitely generated''-assumption and 
correspondingly an ``almost isomorphism''-conclusion for the map induced on homotopy groups. Also there 
is the analogous refinement specifying the range of dimensions where we have to impose the finiteness conditions
in order to obtain the conclusion for a given dimension.

In the case where $G$ is the trivial group assumption (C) is redundant, assumption (D) simply says that the ordinary homology groups
of $X$
are degreewise finitely generated $\bbZ$-modules  and the conclusion is about the map
\[
X_+ \sma \holim_{\cali} \bfE \to \holim_{\cali} X_+ \sma \bfE.
\]

\begin{example}\label{ex:caln}
An important example is the case where the category is~$\caln$, the category whose objects are the non-negative integers
and for which $\mor_{\caln}(i,j)$ is empty if~$i > j$, and consists of one element
if~$i \le j$. In that case a contravariant functor $\bfE: \caln \to \SPECTRA$ is simply an inverse system of spectra
\[
\bfE_0 \leftarrow \bfE_1 \leftarrow \bfE_2\leftarrow \ldots .
\]
There exists a $1$-dimensional model for the contravariant classifying space $E \caln$, compare
Section~\ref{sec: Examples of categories with finite-dimensional classifying space}.
\end{example}

\begin{example} \label{exa: Adams result} 
Illuminating is also the extreme case where the only morphisms in $\cali$ are the identities.
In that case $\cali$ can be identified with an index set $I$ and 
a contravariant functor 
$\cali \to \Omega\text{-}\SPECTRA$ is just a family of
$\Omega$-spectra $\{\bfE_i \mid i \in I\}$ indexed by $I$. 
In the case where $G$ is the trivial group the above specializes to results about the map
\[
X_+ \wedge \prod_{i \in I} \bfE_i ~ \to ~ \prod_{i \in I} \bfR\left(X_+\wedge \bfE_i \right)
\]
and is a variation on~\cite[Theorem~15.2 on page~326]{Adams(1974)}.
\end{example}

Now we would like to briefly explain the concrete situation that motivated the investigations in  this paper.
In~\cite{Lueck-Reich-Rognes-Varisco(2003)}
we study to what extent the methods used by B\"okstedt, Hsiang and Madsen  
in~\cite{Boekstedt-Hsiang-Madsen(1993)} to prove the algebraic $K$-theory Novikov Conjecture can be applied
to prove rational injectivity results for generalized assembly maps in algebraic $K$-theory. 
The main tool in both cases is to 
use the so called cyclotomic trace to topological cyclic homology in order to detect elements in $K$-theory.
The $K$-theoretic Novikov Conjecture says that the ``classical'' assembly map corresponding to the family consisting of the 
trivial subgroup is rationally injective. 
Using generalized assembly maps we can detect a much bigger portion in $K$-theory. The method of comparison via 
the cyclotomic trace 
thus naturally leads to study generalized assembly maps for topological cyclic homology:
\[
EG ( \calf )_+ \sma_{\Or (G) } TC( \bfA ( ? ) ;p) \to 
TC ( \bfA ( G ) ;p).
\]
Here $EG( \calf )$ is the classifying space for a family of
subgroups $\calf$ (compare Section~\ref{sec: The homology of classifying spaces for families}),
$\bfA$~is a fixed coefficient ring(spectrum),
and  for a fixed prime $p$ the $\Or (G)$-spectrum $TC( \bfA( ? );p )$ is defined
as a certain homotopy limit
\[
TC ( \bfA ( ? ) ; p) = \holim_{\calr \calf} THH ( \bfA ( ? ) )^{C_{p^n}}.
\]
The category $\calr \calf$ is described in 
Section~\ref{sec: Examples of categories with finite-dimensional classifying space}
and we verify in Proposition~\ref{the: Ecalc is two-dimensional} 
that it has a $2$-dimensional model for its contravariant classifying space.
Moreover the spectra $THH ( \bfA (G/H))^{C_{p^n}}$ are uniformly bounded below, so that  conditions (A) and (B) 
are satisfied.
In order to study the generalized assembly map it is hence particularly important 
to understand the conditions (C) and (D) 
in the case where $X$ is the classifying space of a family of subgroups.
The following proposition,
which is proven in Section~\ref{sec: The homology of classifying spaces for families},
says that for certain families the conditions can be formulated
entirely in terms of group homology.
\begin{proposition} \label{cor: assumptions of main theorem}
Let $\calf$ be a family of finite subgroups of~$G$ such that the set of conjugacy classes 
$( \calf ) = \{ (H) \; | \; H \in \calf \}$ is finite.

If $H_n(BZ_GH;\bbZ)$ is an almost finitely
generated $\bbZ$-module for all~$n\geq0$ and all~$H \in \calf$,
then the $G$-$CW$-complex $EG(\calf)$ has finite isotropy groups and only
finitely many orbit types, and for any finite group~$H \subset G$
and any~$n\geq0$ the homology group $H_n(Z_GH\backslash
EG ( \calf )^H;\bbZ)$ is an almost finitely generated $\bbZ$-module. 
\end{proposition}

\begin{remark} \label{rem: motivation for almost finitely generated} 
The notions ``almost finitely generated $\bbZ$-module'' and 
``almost isomorphism'' allow to prove 
simultaneously a result about commuting homotopy limits and smash products
as in Addendum~\ref{addendum} and the above 
Proposition~\ref{cor: assumptions of main theorem}. 

If one drops ``almost'', then our main Theorem~\ref{the: commuting wedge and homotopy limit} works, 
but the proof of a corresponding version of Proposition~\ref{cor: assumptions of main theorem} breaks down,
because then $\EGF{G}{\calf}$
may have infinitely many equivariant cells. 
One could get away with the stronger assumption
that $\EGF{G}{\calf}$ has a model with finite skeleta. But that is a lot more restrictive 
and often not true in the applications we have in mind.

On the other hand a version of Proposition~\ref{cor: assumptions of main theorem}
where the assumption and the conclusion is that the rationalized homology groups are finitely generated as $\bbQ$-modules
does exist. But under such an assumption one 
cannot prove an analogue of Theorem~\ref{the: commuting wedge and homotopy limit}, compare Remark~\ref{rem: Q-finitely generated too weak}.
\end{remark}

We finish this introduction with an example that shows that no proper
subset of the conditions appearing in 
Theorem~\ref{the: commuting wedge and homotopy limit} is sufficient.

\begin{example} \label{exa: necessary finiteness conditions on X}
Let $G$ be the trivial group and $\cali$ be the trivial category.
Fix two sequences of integers $0 \le m_0 \le m_1\le m_2 \le \ldots$ and
$0 \le n_0 \le n_1\le n_2 \le \ldots$. Put $X_+ = \bigvee_{i=0}^{\infty} S^{m_i}$ and
$Y_+ = (\coprod_{j=0}^{\infty} S^{n_j})_+$. (Notice that 
here we consider for simplicity pointed spaces instead of spaces with a disjoint base point added).
Let $\bfE$ be an $\Omega$-spectrum. The transformation
$\bft$ appearing in Question~\ref{question} can be identified with the map
\begin{multline*}
\bigvee_{i=0}^{\infty} S^{m_i} \wedge \prod_{j = 0}^{\infty} \map(S^{n_j},\bfE) 
\to 
\prod_{j = 0}^{\infty} \map\left(S^{n_j},\bigvee_{i=0}^{\infty} S^{m_i} \wedge \bfE\right)  
\\
\to 
\prod_{j = 0}^{\infty} \map\left(S^{n_j};
\bfR\left(\bigvee_{i=0}^{\infty} S^{m_i} \wedge \bfE\right)\right).
\end{multline*} 
Recall that for a spectrum $\bfF$ we have $\pi_p(S^n \wedge \bfF) = \pi_{p-n}(\bfF)$ and  
$\pi_p(\map(S^n,\bfF)) = \pi_{p+n}(\bfF)$, and  for a collection
$\{X_i\mid i \in I\}$ of pointed $CW$-complexes the canonical map
$$\bigoplus_{i \in I} \pi_p\left(X_i \wedge \bfF\right) 
~ \xrightarrow{\cong} ~ 
\pi_p\left(\bigvee_{i \in I} X_i \wedge \bfF\right)$$
is bijective.
Using Lemma~\ref{lem: commuting pi and products for Omega-spectra} and
Lemma~\ref{lem: map respects Omega-spectra} one can identify the map above with
the canonical map
$$\bigoplus_{i=0}^{\infty} \prod_{j=0}^{\infty} \pi_{n_j-m_i+p}(\bfE) 
~ \to ~
 \prod_{j=0}^{\infty} \bigoplus_{i=0}^{\infty} \pi_{n_j-m_i+p}(\bfE).$$
This map is always injective. But it is surjective if and only if there exists
$i_0$ such that for all $i \ge i_0$ and all $j$ we have
$\pi_{n_j-m_i+p}(\bfE) = 0$. Note that condition~(A) corresponds to the existence of a constant $D$ such that
$n_j \leq D$ for all $j$. 
Condition~(B) says that there exists a constant $N$ such that $\pi_q( \bfE ) = 0$ for $q \leq N$. 
Condition~(C) is redundant in the non-equivariant case and
condition~(D) is satisfied if and only if $\lim_{i \to \infty} m_i = \infty$. 
It is then easy to find examples showing that we cannot drop any of the assumptions (A), (B) or (D).
\end{example}

\typeout{--------------------  Basics about spectra  --------------------}

\section{Basics about spectra}
\label{sec: Basics about spectra}

The question we are studying is not even well posed in the stable homotopy category of spectra.
We are hence forced to explain what we mean by a spectrum, and in what sense we understand the 
basic constructions concerning spectra.

Throughout this paper we work in the category of compactly generated
spaces (see~\cite{Steenrod(1967)}, \cite[I.4]{Whitehead(1978)}). 
So space means compactly generated space and all constructions like
mapping spaces and products are to be understood in this category. 
We will always assume that the inclusion of the base point into a pointed space is a cofibration
and that maps between pointed spaces preserve the base point.

We define the category $\SPECTRA$
as follows. A {\em spectrum}
$\bfE = \{(E(n),\sigma(n)) \mid n = 0, 1, 2, \dots \}$
is a sequence of pointed spaces
$\{E(n) \mid n =0 ,1 , 2, \dots\}$ together with pointed maps
(called {\em structure maps})
$\sigma(n) : E(n) \wedge S^1 \to E(n+1)$.
A \emph{map} of spectra (sometimes also called function in the literature)
$\bff\colon  \bfE \to \bfE^{\prime}$ is a sequence of
maps of pointed spaces
$f (n)\colon E(n) \to E^{\prime}(n)$
that are compatible with the structure maps $\sigma(n)$, i.e.~we have
$f(n+1) \circ \sigma(n) ~ = ~ \sigma^{\prime}(n) \circ \left(f(n) \wedge \id_{S^1}\right)$
for all~$n$.
This should not be confused with the notion of map of spectra in
the stable homotopy category (see~\cite[III.2]{Adams(1974)}).
Recall that the homotopy groups of a spectrum
are defined as
\[
\pi_p(\bfE) = \colim_{k \to \infty} \pi_{p+k}(E(k)),
\]
where the maps in this system are given by the composition
$$\pi_{p+k}(E(k)) \longrightarrow \pi_{p+k+1}(E(k)\wedge S^1)
\xrightarrow{\sigma (k)_*} \pi_{p+k+1}(E(k +1))$$
of the suspension homomorphism and the
homomorphism induced by the structure map.
A {\em weak equivalence}
of spectra is a map $\bff \colon \bfE \to \bfF$ 
inducing an isomorphism on all homotopy
groups.  A spectrum $\bfE$ is called {\em $\Omega$-spectrum}
if the adjoint of each structure  map 
$$\overline\sigma(n)\colon E(n) \to \Omega E(n+1) = \map (S^1,E(n+1))$$
is a weak homotopy equivalence of spaces. We denote by
$\Omega$-$\SPECTRA$ the corresponding full subcategory of $\SPECTRA$.

In order to define the covariant \emph{fibrant replacement} functor
\[
  \bfR \colon \SPECTRA \to \Omega\text{-}\SPECTRA
\]
we first need to recall the definition and a key property of sequential
homotopy colimits.  Consider a sequence of pointed spaces and maps
$
  X_0 \to X_1 \to X_2 \to \ldots
$
This can be thought of equivalently as a covariant functor~$X$ from
the category~$\caln$ to pointed spaces (see Example~\ref{ex:caln}).
Its homotopy colimit is defined as the pointed space
\[
  \hocolim_{n \to \infty} X_n = \hocolim_\caln X =
  E\caln_+ \sma_\caln X ,
\]
compare Definition~\ref{def:holim}.  Using the model for~$E\caln$
described in Example~\ref{ex:Ecaln} one gets the well-known reduced
mapping telescope construction.  Since any map from a compact space to
$\hocolim_{n \to \infty} X_n$ factors for some~$N$ through the natural
map $X_N \to \hocolim_{n \to \infty} X_n$, we have that for any~$k\geq0$
\begin{equation}\label{eq:pihocolim}
  \pi_k(\hocolim_{n \to \infty} X_n)
  \cong
  \colim_{n \to \infty} \pi_k(X_n) .
\end{equation}

Now for any spectrum~$\bfE$ we define its fibrant replacement as
the spectrum $\bfR \bfE$ whose $k$-th pointed space is
$\hocolim_{n \to \infty} \Omega^nE(n+k)$.
Here the homotopy colimit is taken with respect to the system of maps
\[
  \Omega^n\overline\sigma(n+k)\colon  \Omega^nE(n+k) \to 
  \Omega^n \Omega E(n+k+1)
  = \Omega^{n+1}E((n+1)+k)
\]
where $\overline\sigma(k) :E(k) \to \Omega E(k+1)$ is the adjoint of
the structure map of $\bfE$.  The adjoints of the structure maps of
$\bfR\bfE$ are given by the composition
\begin{multline*}
\hocolim_{n \to \infty} \Omega^nE(n+k) 
\xrightarrow{\hocolim_{n \to \infty} \Omega^{n}\overline\sigma(k+n)}
\hocolim_{n \to \infty} \Omega^n\Omega E(n+k+1) 
\\
=
\hocolim_{n \to \infty} \Omega \Omega^n E(n+(k+1)) \xrightarrow{\mu}
\Omega \hocolim_{n \to \infty} \Omega^n E(n+(k+1)),
\end{multline*}
where the first map is the obvious weak equivalence induced by shifting
the system, and $\mu$ is the canonical map, which is also a weak
equivalence as one sees applying~\eqref{eq:pihocolim} above.  Hence
$\bfR\bfE$ is indeed an $\Omega$-spectrum.
There is a natural map
\begin{eqnarray}
& \bfr(\bfE) \colon  \bfE \to \bfR \bfE &
\label{bfr}
\end{eqnarray}
coming from the canonical maps
$E(k) \to \hocolim_{n \to \infty} \Omega^nE(n+k)$.
Using~\eqref{eq:pihocolim} one computes that $\bfr(\bfE)$ is a weak
equivalence.

Many constructions for spaces can be carried over to spectra by applying them levelwise,
i.e.~to each of the individual spaces $E(n)$. Usually there is then an obvious way to define
either the structure maps or their adjoints. For example the \emph{product}
$\prod_{i \in I} \bfE_i$ of a 
collection of spectra $\{\bfE_i \mid i \in I\}$
is the spectrum whose $n$-th
space is $\prod_{i \in I} E_i(n)$. For the structure maps observe that $\Omega$ commutes with 
products.
For us the main reason to consider 
$\Omega$-spectra is the following lemma, which will play an important role below and
does not hold for  arbitrary spectra.

\begin{lemma} \label{lem: commuting pi and products for Omega-spectra}
Let $\{\bfE_i \mid i \in I\}$ be a family of $\Omega$-spectra, where
$I$ is an arbitrary index set.
Then the canonical map induced by the various projections
$\pr_i\colon  \prod_{i \in I} \bfE_i \to \bfE_i$ 
\[
  \pi_n\left(\prod_{i \in I} \bfE_i\right) \xrightarrow{\cong}
  \prod_{i \in I} \pi_n(\bfE_i)
\]
is bijective for all~$n \in \bbZ$.
\end{lemma}

\begin{proof}
Notice that for an $\Omega$-spectrum $\bfF$ the canonical map
$$\pi_n(F(0)) \to \pi_n(\bfF) = \colim_{k \to \infty} \pi_n(\Omega^{k}F(n+k))$$
is bijective and that $\prod_{i \in I} \bfE_i$ is an $\Omega$-spectrum.
Hence the claim follows from the corresponding statement for spaces.
\end{proof}

The following example shows that
Lemma~\ref{lem: commuting pi and products for Omega-spectra}
does not hold without the assumption that each
$\bfE_i$ is an $\Omega$-spectrum.

\begin{example} \label{exa: assumption Omega-spectrum is necessary}
Let $\bfE$ be an $\Omega$-spectrum such that $A=\pi_0 ( \bfE ) \cong \pi_k ( E (k) )$ is a 
nontrivial abelian group. Denote by $\bfE_{(i, \infty)}$ the spectrum obtained from $\bfE$ by replacing the 
spaces $E (0)$, $E (1)$, $\dots$, $E(i)$ by the one-point space $\{ \pt \}$.
We have 
\[
\pi_k ( E_{(i, \infty)}(k) ) = \left\{ \begin{array}{cl} 0 & \mbox{ if } k \leq i \\ A & \mbox{ if } k > i \end{array} \right.
\]
and the maps $\pi_k ( E_{(i, \infty)} (k) ) \to \pi_{k+1} ( E_{(i, \infty)} (k+1) )$ are isomorphisms
for $k > i$. We see that for all $i$ we have $\pi_0 ( \bfE_{(i , \infty)} ) = A$ and 
\[
\pi_0 ( \prod_{i \in \bbN} \bfE_{(i , \infty)} ) = \colim_k \pi_k ( \prod_{i \in \bbN} E_{(i, \infty)} (k) )
= \colim_k \bigoplus_{i=1}^k A = \bigoplus_{i=1}^{\infty} A.
\]
The natural map
\[
\pi_0 ( \prod_{i \in \bbN} \bfE_{(i, \infty)} ) \to \prod_{i \in \bbN} \pi_0 ( \bfE_{(i,\infty)} )
\]
can be identified with the natural inclusion 
$\bigoplus_{i \in \bbN} A \to \prod_{i \in \bbN} A$
and is not an isomorphism.
\end{example}

Now we fix some notations and conventions about homotopy pullbacks and
homotopy pushouts of spectra. Consider a commutative square $D_{\bfE}$ of spectra
\[
\xymatrix{
  \bfE_0 \ar[r]^{\bff_1} \ar[d]_{\bff_2} & \bfE_1    \ar[d]^{\bfg_1} \\
  \bfE_2 \ar[r]_{\bfg_2}                 & \bfE_{12}.
}
\]
We denote by $\hopb{\bfE_2}{\bfE_{12}}{\bfE_1}$
the levelwise homotopy pullback spectrum
and by $\hopo{\bfE_2}{\bfE_0}{\bfE_1}$
the levelwise homotopy pushout spectrum,
i.e.~the $k$-th spaces are given by 
the homotopy pullback respectively the homotopy pushout
of the corresponding diagrams of pointed spaces.
For the structure maps use the fact that homotopy pullbacks commute with~$\Omega$ and
homotopy pushouts commute with~$S^1 \wedge -$ up to 
natural homeomorphisms.
There are canonical maps of spectra
\[
\bfa \colon \bfE_0 \to \hopb{\bfE_2}{\bfE_{12}}{\bfE_1} , \qquad
\bfb \colon \hopo{\bfE_2}{\bfE_0}{\bfE_1} \to \bfE_{12} .
\]
We call $D_\bfE$ \emph{homotopy cartesian} if the map~$\bfa$ is a
weak equivalence of spectra; dually we call $D_\bfE$
\emph{homotopy cocartesian} if~$\bfb$ is a weak equivalence. We use the analogous terminology
for spaces.

\begin{remark} \label{rem: homotopy pullback with about R}
A homotopy pullback is a special case of a homotopy limit (the indexing category is 
$2\!\ot\!12\!\to\!1$ for contravariant functors).
In Definition~\ref{def:holim} below, we discuss 
that in order to define a homotopy limit of a diagram of spectra, one
should first replace all spectra in sight with $\Omega$-spectra by
applying the fibrant replacement functor~$\bfR$, and then take the
levelwise homotopy limit.  At first glance this seems to be
inconsistent with the definition of homotopy pullback given above.
But Lemma~\ref{lem: hopushout and hopullback for spectra} below implies in
particular that the canonical map~$\bfr$ induces a weak equivalence of
spectra
\[\hopb{\bfE_2}{\bfE_{12}}{\bfE_1} \to 
\hopb{\bfR\bfE_2}{\bfR\bfE_{12}}{\bfR\bfE_1}.\]
\end{remark}

Finally we recall the well-known fact that a commutative square of
spectra is homotopy cartesian if and only if it is homotopy
cocartesian.  
We would like to stress that it is important for later proofs 
to have such a statement for $\hopb{\bfE_2}{\bfE_{12}}{\bfE_1}$ 
as opposed to $\hopb{\bfR\bfE_2}{\bfR\bfE_{12}}{\bfR\bfE_1}$.
We have not been able to find an explicit reference
(in particular none dealing with this extra subtlety)
in the literature, therefore we include a proof.

\begin{lemma} \label{lem: hopushout and hopullback for spectra}
A commutative square~$D_{\bfE}$ of spectra is homotopy cocartesian if
and only if it is homotopy cartesian. In this case there is a natural
long exact Mayer-Vietoris sequence
\begin{multline*}
\ldots \xrightarrow{\partial_{n+1}} \pi_n(\bfE_0) 
\xrightarrow{\pi_n(\bff_1) \oplus \pi_n(\bff_2)}
\pi_n(\bfE_1) \oplus \pi_n(\bfE_2) 
\\
\xrightarrow{\pi_n(\bfg_1) - \pi_n(\bfg_2)} \pi_n(\bfE_{12}) 
\xrightarrow{\partial_n}
\pi_{n-1}(\bfE_0) \xrightarrow{\pi_{n-1}(\bff_1) \oplus
  \pi_{n-1}(\bff_2)}
\ldots
\end{multline*}
\end{lemma}

\begin{example} \label{exa: wedge and product}
A special case of this Lemma yields the following well-known fact:
For any pair of spectra~$\bfE_1$ and~$\bfE_2$ the natural map
$\bfE_1 \vee \bfE_2 \to \bfE_1 \times \bfE_2$
is a weak equivalence.
\end{example}

\begin{proof}[Proof of Lemma~\ref{lem: hopushout and hopullback for spectra}]
The key step in the proof is to show the following
\begin{itemize}
\item[\emph{Claim:}] \emph{If a commutative square of spectra is
    levelwise a homotopy cocartesian square of spaces, then it is
    homotopy cartesian.}
\end{itemize}
We first explain why this claim implies the lemma.
Any levelwise homotopy cartesian square of spectra yields a natural
long exact Mayer-Vietoris sequence.  Using the natural map~$\bfa$ we
get a natural Mayer-Vietoris sequence for every homotopy cartesian
square.  For every commutative square~$D_\bfE$ we consider the levelwise
homotopy cocartesian square~$D_\text{po}$
\[
\xymatrix{
  \bfE_0 \ar[r]^-{\bff_1} \ar[d]_{\bff_2} & \bfE_1 \ar[d]^{\overline\bfg_1} \\
  \bfE_2 \ar[r]_-{\overline\bfg_2}        & \hopo{\bfE_2}{\bfE_0}{\bfE_1} .
}
\]
The claim implies that $D_\text{po}$ is homotopy cartesian, and
therefore there is a long exact Mayer-Vietoris sequence associated
to~$D_\text{po}$.

Now assume that~$D_\bfE$ is homotopy cartesian.  In this case both
$D_\text{po}$ and~$D_\bfE$ induce long exact Mayer-Vietoris sequences,
and the map~$\bfb$ induces a natural transformation between them.  The
Five-Lemma then implies that $\bfb$ is a weak equivalence of spectra,
i.e.~$D_\bfE$ is homotopy cocartesian.  Conversely, assume that
$D_\bfE$ is homotopy cocartesian.  Using the natural weak
equivalence~$\bfb$ we get a long exact Mayer-Vietoris sequence
for~$D_\bfE$ from the one for~$D_\text{po}$.  Then a Five-Lemma
argument shows that the map~$\bfa$ is a weak equivalence,
i.e.~$D_\bfE$ is homotopy cartesian.

It remains to prove the claim.  This is based on the following
well-known stable range statement about homotopy cocartesian squares
of spaces (a version of Blakers-Massey homotopy excision theorem, see
for example~\cite[Section~2, in particular
Theorem~2.3]{Goodwillie(1991)}).  Assume that the
commutative square of spaces
\[
\xymatrix{
  X_0 \ar[r]^{f_1} \ar[d]_{f_2} & X_1 \ar[d]^{g_2} \\
  X_2 \ar[r]_{g_1} & X_{12}
}
\]
is homotopy cocartesian,
and that $f_1$~is $k_1$-connected and $f_2$~is $k_2$-connected.
Then the natural map :
$a\colon X_0 \to \hopb{X_2}{X_{12}}{X_1}$
is at least ($k_1+k_2-1$)-connected.

Now suppose we are given a levelwise homotopy cocartesian
square~$D_\bfE$ of spectra.  Since homotopy
pushouts commute with suspensions we have for every $n,k\in\bbN$ a
homotopy cocartesian square of spaces
\[
\xymatrix{
  \Sigma^k E_0(n) \ar[r]^{\Sigma^k f_1(n)} \ar[d]_{\Sigma^k f_2(n)} 
  & \Sigma^k E_1(n) \ar[d]^{\Sigma^k g_2(n)} \\
  \Sigma^k E_2(n) \ar[r]_{\Sigma^k g_1(n)} & \Sigma^k E_{12}(n) .
}
\]
Now consider the following diagram
  \[
  \xymatrix@R=3.5ex{
  \Sigma^k E_0(n) \ar@{=}[d] \ar[r]^-{\Sigma^k a(n)} & 
  \Sigma^k \hopb{E_2(n)}{E_{12}(n)}{E_1(n)} \ar[d]^\eta \\
  \Sigma^k E_0(n) \ar[d]_{\sigma} \ar[r]^-{a_{k,n}} &
  \hopb{\Sigma^k E_2(n)}{\Sigma^k E_{12}(n)}{\Sigma^k E_1(n)} \ar[d]^\mu \\
  E_0(n+k) \ar[r]^-{a(n+k)} & \hopb{E_2(n+k)}{E_{12}(n+k)}{E_1(n+k)}
  }
\]
where $\sigma$ comes from the structure maps of~$\bfE_0$, $\eta$~is
the natural map, and~$\mu$~the map induced by the structure maps
of~$\bfE_2$, $\bfE_{12}$ and $\bfE_1$.  One checks that the diagram
commutes and that the composition~$\mu\circ\eta$ is the structure map
of the spectrum~$\hopb{\bfE_2}{\bfE_{12}}{\bfE_1}$ (recall that for
the homotopy pullback spectrum one actually defines the
adjoints of the structure maps, using the fact that homotopy pullbacks
and loops commute up to natural homeomorphism).  From the
Blakers-Massey connectivity statement recalled above we get that the
natural map~$a_{k,n}$ in the middle row of the diagram is at least
$(2k-3)$-connected.
Hence we can apply Lemma~\ref{lem:connectivity} below to conclude that
the natural map of spectra
$\bfa\colon\bfE_0\to\hopb{\bfE_2}{\bfE_{12}}{\bfE_1}$
is a weak equivalence, i.e.~$D_\bfE$~is homotopy cartesian.
\end{proof}

\begin{lemma}\label{lem:connectivity}
  Let~$\bff\colon\bfE\to\bfF$ be a map of spectra
  and $l\in\bbZ$~be an integer.
  Suppose that for all $n\in\bbN$ there exist a~$k=k(l,n)\in\bbN$
  and a commutative diagram
  \[
  \xymatrix@R=2ex{
    \Sigma^k E(n) \ar[d] \ar@/_4.5ex/[dd] \ar[r]^-{\Sigma^k f(n)} &
    \Sigma^k F(n) \ar[d] \ar@/^4.5ex/[dd] \\
    E_{k,n} \ar[d] \ar[r]^-{f_{k,n}} &
    F_{k,n} \ar[d] \\
    E(k+n) \ar[r]^-{f(k+n)} & F(k+n)
  }
  \]
  where $f_{k,n}$ is $(l+k+n)$-connected, and the outer vertical maps
  are given by the structure maps of the spectra.
  Then~$\bff\colon\bfE\to\bfF$ is $l$-connected.
\end{lemma}

We omit the elementary proof, which just uses the definition of
homotopy groups of spectra, basic properties of sequential colimits,
and diagram chases.

\typeout{--------------------  Basics about ... over a category  --------------------}

\section{Basics about spaces, spectra and modules over a category}
\label{sec: Basics about modules and spectra over a category}

In this section we recall some basic notions about spaces, spectra and
modules over a small category~$\cali$.  More information for spaces and
spectra can be found in~\cite{Davis-Lueck(1998)} and for modules
in~\cite[II.9]{Lueck(1989)}.

A \emph{covariant (contravariant) pointed $\cali$-space $X$}
is a covariant (contravariant) functor
from $\cali$ to the category of pointed (compactly generated) spaces.
A map between $\cali$-spaces is a natural transformation
of such functors.  For pointed $\cali$-spaces $X$ and $Y$ of the
same variance let $\map_{\cali}(X,Y)$
be the space of pointed maps from $X$
to $Y$ with the subspace topology coming from the obvious inclusion
into $\prod_{c \in \obj(\cali)} \map(X(c),Y(c))$. 
Let $X$ be a contravariant and~$Y$ a covariant $\cali$-space. Define
their \emph{smash product} to be the space
\begin{eqnarray}
X \sma_{\cali} Y ~ = ~  \bigvee_{c \in \obj(\cali)} X(c) \wedge Y(c)/\sim,
\label{X otimes_{cali} Y}
\end{eqnarray}
where $\sim$ is the equivalence relation generated by
$(x\phi,y) \sim (x,\phi y)$
for all morphisms  $\phi\colon c \to d$ in $\cali$ and all points
$x \in X(d)$ and $y \in Y(c)$. Here $x\phi$ stands for $X(\phi)(x)$ and $\phi y$
for $Y(\phi )(y)$. If $\cali$ is the trivial category with only one object
and only one morphism, then  $X$ and $Y$ are just pointed spaces and 
$X \sma_{\cali} Y$ is the ordinary smash product $X \wedge Y$.
For an ordinary pointed space $Z$ let $\map(Y,Z)$ denote the
obvious contravariant pointed $\cali$-space whose value at an object~$c$ is  the
mapping space $\map(Y(c),Z)$ of base point preserving maps. 
One easily checks the following lemma.

\begin{lemma} \label{lem: adjointness of tensor product}
Let $X$ be a contravariant pointed  $\cali$-space, $Y$ be a covariant
pointed $\cali$-space and $Z$ be a pointed space. 
Then there is a homeomorphism
$$T \colon  \map(X \sma_{\cali} Y,Z) ~ \xrightarrow{\cong} ~
\map_{\cali}(X,\map(Y,Z))$$
which is natural in $X$, $Y$ and $Z$.
\end{lemma}

A \emph{covariant (contravariant) $\cali$-spectrum} is a covariant (contravariant) functor
$\bfE\colon \cali \to \SPECTRA$. 
A  \emph{covariant (contravariant) $\cali$-$\Omega$-spectrum} is a covariant (contravariant) functor
$\bfE\colon \cali \to \Omega\text{-}\SPECTRA$. Given a contravariant pointed
$\cali$-space $X$, a covariant $\cali$-spectrum $\bfE$ and a contravariant $\cali$-spectrum $\bfF$ we 
obtain spectra 
\[
X \sma_{\cali} \bfE \quad \mbox{ and } \quad \map_{\cali} ( X , \bfF )
\]
in the obvious way by applying the constructions discussed above 
levelwise. Notice that for an $\cali$-$\Omega$-spectrum $\bfE$ the spectrum
$X \sma_{\cali} \bfE$ is not necessarily an $\Omega$-spectrum.

\begin{definition}
\label{def: cali-CW-complexes}
A {\em contravariant $\cali$-$CW$-complex}
$X$ is a contravariant unpointed
$\cali$-space $X \colon \cali \to \SPACES$ together with a filtration
$$\emptyset = X_{-1} \subset X_0 \subset X_1 \subset X_2  \subset\ldots
\subset X_n \subset \ldots \subset X = \bigcup_{n \ge 0} X_n$$
such that $X = \colim_{n \to \infty} X_n$ and for any $n \ge 0$ the
{\em $n$-skeleton} $X_n$  is obtained from the
$(n-1)$-skeleton $X_{n-1}$ by attaching free contravariant
$\cali$-$n$-cells,
i.e.~there exists a pushout of $\cali$-spaces of the form
\[
\xymatrix{
  \coprod_{i \in I_n} \mor_{\cali}(-,c_i) \times S^{n-1} \ar[d] \ar[r]
    & X_{n-1} \ar[d] \\
  \coprod_{i \in I_n} \mor_{\cali}(-,c_i) \times D^n            \ar[r]
    & X_n ,
}
\]
where the vertical maps are inclusions, $I_n$ is an index set, and
the $c_i$ are objects of~$\cali$.
(In~\cite[Definition~3.2 on page~221]{Davis-Lueck(1998)} 
this is called a free $\cali$-$CW$-complex, but we will drop the word free here.)
\end{definition}

\begin{definition} \label{def: Ecalc} 
A contravariant $\cali$-$CW$-complex $E\cali$ whose 
evaluation at each object of $\cali$ is a contractible space is called a
\emph{contravariant classifying space} for $\cali$.
\end{definition}

\begin{remark}
Any two models for~$E\cali$ are $\cali$-homotopy equivalent. 
A bar-construc\-tion model for $E\cali$, which is functorial in~$\cali$,
is given in~\cite[page~230]{Davis-Lueck(1998)}.
\end{remark}

\begin{definition}\label{def:holim}
If $X$ is a covariant pointed $\cali$-space and $Y$~a contravariant pointed $\cali$-space
the \emph{homotopy colimit} respectively the \emph{homotopy limit} are defined as the pointed spaces given by
\[
\hocolim_{\cali} X =  E\cali_+ \sma_{\cali} X \quad \mbox{ respectively } \quad \holim_{\cali} Y =\map_{\cali}(E\cali_+,Y),
\]
compare~\cite[Definition~3.13 on page~225]{Davis-Lueck(1998)}. For a covariant $\cali$-spectrum $\bfE$ and a contravariant $\cali$-spectrum
$\bfF$ we define similarly
\[
\hocolim_{\cali} \bfE  =  E\cali_+ \sma_{\cali} \bfE 
\quad \mbox{ and  } \quad 
\holim_{\cali} \bfF =\map_{\cali}( E\cali_+,\bfR \circ \bfF).
\]
But note the replacement functor that appears in the homotopy limit and the fact that even if
$\bfE$ is an $\cali$-$\Omega$-spectrum the homotopy colimit need not be an $\Omega$-spectrum.
\end{definition}

\begin{remark}\label{rem:propholim}
The homotopy type of the homotopy colimit or homotopy limit does not
depend on the choice of the model for~$E\cali$.
If the bar-construction model for~$E\cali$ is used, then our
definition agrees with the obvious modification
of the one given in~\cite[Chapters~XI and~XII]{Bousfield-Kan(1972)}
for diagrams of simplicial sets.
Some properties of homotopy colimits and homotopy limits become
obvious with our definition.
For instance $\Omega$ commutes with homotopy limits of spaces because
of the adjunctions
\[
  \map(S^1,\map_{\cali}(E\cali_+,X))
  \cong
  \map_{\cali}(S^1\sma_{\cali}E\cali_+,X)
  \cong
  \map_{\cali}(E\cali_+,\map(S^1,X))
\]
(see~\cite[Lemma~1.6 on page~206]{Davis-Lueck(1998)}).
This fact is used to define the structure maps of the homotopy limit of spectra.
\end{remark}

We collect some results that will be needed below
(compare~\cite[Theorem~3.11 on page~224 and Lemma~4.6 on page~234]{Davis-Lueck(1998)}).

\begin{lemma} \label{lem: weak equivalences and compatibility with otimes and map}
Let $Y$ be a contravariant $\cali$-$CW$-complex.

\begin{enumerate}

\item \label{lem: weak equivalences and compatibility with otimes and map: spaces}

Let $f\colon Z_1 \to Z_2$ be a map of covariant or contravariant
respectively $\cali$-spaces such 
that $f(c)$ is a weak equivalence for all objects $c$ in $\cali$.
Then the induced maps
$$Y \sma_{\cali} Z_1 \to Y \sma_{\cali} Z_2$$
and
$$\map_{\cali}(Y,Z_1) \to \map_{\cali}(Y,Z_2)$$
respectively are weak equivalences.

\item \label{lem: weak equivalences and compatibility with otimes and map: spectra and wedge}

Let $f\colon \bfE \to \bfF $ be a map of covariant  $\cali$-spectra such 
that $f(c)$ is a weak equivalence for all objects $c$ in $\cali$. 
Then the induced map
$$Y \sma_{\cali} \bfE\to Y \sma_{\cali} \bfF$$
is a weak equivalence.

\item \label{lem: weak equivalences and compatibility with otimes and map: spectra and map}

Let $f\colon \bfE \to \bfF $ be a map of contravariant  $\cali$-$\Omega$-spectra such 
that $f(c)$ is a weak equivalence for all objects $c$ in $\cali$.
Then the induced map
$$\map_{\cali}(Y,\bfE) \to \map_{\cali}(Y,\bfF)$$
is a weak equivalence.

\end{enumerate}

\end{lemma}

\begin{lemma}
\label{lem: map respects Omega-spectra}
Let $\bfE\colon  \cali \to \Omega\text{-}\SPECTRA$ be a contravariant functor and
let $Y$ be a contravariant $\cali$-space. Then
$\map_{\cali}(Y_+,\bfE)$ is an $\Omega$-spectrum.
\end{lemma}

\begin{proof}
The adjoint of the $k$-th structure map of $\map_{\cali}(Y_+,\bfE)$ is the composition
$$\map_{\cali}(Y_+,E(k)) \xrightarrow{\map_{\cali}(\id,\overline\sigma(k))} 
\map_{\cali}(Y_+,\Omega E(k+1)) \xrightarrow{\cong} \Omega
\map_{\cali}(Y_+,E(k+1)),$$
where the first map is a weak homotopy equivalence by 
Lemma~\ref{lem: weak equivalences and compatibility with otimes and map}~(i)
and the second one is the canonical homeomorphism
(compare Remark~\ref{rem:propholim}).
\end{proof}

The constructions discussed above for spaces and spectra all have analogues in the world
of modules and chain complexes. We will usually refer to this analogy as the ``linear case''.
We collect some basics. Fix a commutative associative ring $R$ with
unit. A \emph{covariant (contravariant) $R\cali$-module} is a covariant (contravariant)
functor from  $\cali$ to the category of $R$-modules.
Morphisms are natural transformations. 

Given a contravariant $R\cali$-module $M$ and a 
covariant $R\cali$-module $N$, their \emph{tensor product} is the $R$-module
\[
M \otimes_{R \cali} N ~ = ~  \bigoplus_{c \in \obj(\cali)} M(c) \otimes_R N(c)/\sim,
\]
where $\sim$ is defined similar as in~\eqref{X otimes_{cali} Y}.
Given two $R\cali$-modules $M$ and $N$ of the same variance, define
$\hom_{R\cali}(M,N)$ to be the $R$-module of natural transformations from $M$ to $N$.
The obvious analogue of  Lemma~\ref{lem: adjointness of tensor product} is true for
$\otimes_{R\cali}$ and  $\hom_{R\cali}$.

The category $R\cali \text{-}\MODULES$ of covariant (contravariant) $R\cali$-modules inherits 
the structure of an abelian category from the category of $R$-modules. In particular
the notion of a projective $R\cali$-module makes sense. 
There are enough injectives and projectives in $R\cali
\text{-}\MODULES$ (see~\cite[9.16 on page~167]{Lueck(1989)},
\cite[Exercise~2.3.5 on page~42 and Example~2.3.13 on page~43]{Weibel(1994)}). Hence 
standard homological algebra applies to $R\cali\text{-}\MODULES$. For instance,
one can define the $R$-module $\Tor^{R\cali}_p(M,N)$ for a  contravariant $R\cali$-module $M$ and a 
covariant $R\cali$-module $N$.

An \emph{$R\cali$-chain complex} is a chain complex in 
$R\cali\text{-}\MODULES$, or, equivalently, a covariant (contravariant) functor from
$\cali$ to the category of $R$-chain complexes.
We denote by~$\ch(R\cali\text{-}\MODULES)$ the category of $R\cali$-chain complexes.

We now specialize to the case where~$\cali$ is the \emph{orbit category}~$\Or(G)$,
whose objects are the homogeneous
$G$-spaces $G/H$ and whose morphisms are $G$-maps. A (left) $G$-space $X$ defines a
contravariant $\Or (G)$-space 
$$X \colon \Or (G) \to \SPACES, \hspace{5mm} G/H \mapsto \map_G(G/H,X) = X^H.$$
If $X$ is a $G$-$CW$-complex, then $X$ considered as a functor is a contravariant $\Or(G)$-$CW$-complex
(see~\cite[Theorem~7.4 on page~250]{Davis-Lueck(1998)}). 
Every $\Or(G)$-$CW$-complex defines a functor from~$\Or(G)$ to the category
of $CW$-complexes, which we can compose with the functor cellular chain complex.
Thus we associate to an $\Or(G)$-$CW$-complex $X$ its \emph{cellular $\bbZ\Or(G)$-chain complex} 
$C_*^{\bbZ\Or(G)}(X)$. This is always a free $\bbZ\Or(G)$-chain complex in the sense
of~\cite[9.16 on page~167]{Lueck(1989)} and in particular projective.

Furthermore we will also need the category $\Sub(G)$, which is a quotient category of $\Or ( G )$.
Its objects are the subgroups
$H$ of $G$. For two subgroups $H$ and
$K$ of $G$ denote by $\conhom_G(H,K)$ the set
of group homomorphisms $f\colon H \to K$,
for which there exists an element $g \in G$
with $gHg^{-1} \subset K$  such that
$f$ is given by conjugation with $g$, i.e.
$f = c(g)\colon H \to K, \hspace{3mm} h \mapsto ghg^{-1}$.
Notice that $c(g) = c(g')$ holds for two elements $g,g' \in G$ with
$gHg^{-1} \subset K$ and $g'H(g')^{-1} \subset K$ if and only if
$g^{-1}g'$ lies in the centralizer
$Z_GH = \{g \in G \mid gh=hg \mbox{ for all } h \in H\}$
of $H$ in $G$. The group of inner automorphisms of $K$ acts on
$\conhom_G(H,K)$ from the left by composition. Define the set of morphisms
$$\mor_{\Sub(G)}(H,K) ~ = ~ \Inn(K)\backslash \conhom_G(H,K).$$

There is a natural projection 
\begin{equation}
\pr\colon  \Or(G) \to \Sub(G)
\label{pr colon Or(G) to Sub(G)}
\end{equation}
which sends a homogeneous space $G/H$ to $H$.
Given a $G$-map $f\colon G/H \to G/K$, we can choose
an element $g \in G$ with $gHg^{-1} \subset K$
and $f(g'H) = g'g^{-1}K$. Then
$\pr(f)$ is represented by $c(g)\colon H \to K$. Notice that
$\mor_{\Sub(G)}(H,K)$ can be identified with the quotient
$\mor_{\Or(G)}(G/H,G/K)/Z_GH$,
where $g \in Z_GH$ acts on $\mor_{\Or(G)}(G/H,G/K)$
by composition with
$R_{g^{-1}}\colon G/H \to G/H, \hspace{3mm} g'H \mapsto g'g^{-1}H$.

Recall that a functor $\bfF \colon \GROUPOIDS \to \SPECTRA$ 
is a homotopy functor if it sends an equivalence of groupoids to a weak equivalence of spectra.
This is equivalent to requiring that naturally isomorphic functors between groupoids induce
the same map after applying $\pi_{\ast} ( \bfF (-))$.
Recall also the transport groupoid functor~$\calg^G$ defined before
Theorem~\ref{the: commuting wedge and homotopy limit}.

\begin{lemma} \label{lemma-factorize-over-groupoids}
Let $\bfF \colon \GROUPOIDS \to \SPECTRA$ be a homotopy functor, then
the functor 
$\pi_{n} ( \bfF \circ \calg^G ( - ) )$
factorizes over the projection $\pr \colon \Or (G) \to \Sub (G)$.
\end{lemma}
\begin{proof}
It suffices to show that the identity and $\calg^G ( R_{d^{-1}} ) \colon \calg^G ( G/H) \to \calg^G ( G/H)$
are naturally isomorphic for $d \in Z_G H$.
A natural transformation is given at the object $gH$ by $gd^{-1}g^{-1} \colon gH \to gd^{-1} H$.
\end{proof}

\typeout{--------------------  Almost trivial and almost f.g.  --------------------}

\section{Almost trivial and almost finitely generated modules}
\label{sec: Almost trivial and almost finitely generated modules}

In this section we discuss some elementary properties of the
notions ``almost finitely generated modules'' and ``almost
isomorphisms'' that were used in Addendum~\ref{addendum}.
The main aim is to verify that we have the fundamental properties that
are necessary to do spectral sequence comparison arguments.

Let $R$ be an integral domain, i.e.~an 
associative commutative ring with unit such that $R$ has no non-trivial
zero-divisor. An $R$-module $M$ is
called \emph{faithful} if its annihilator ideal 
$\Ann_R(M) = \{r \in R \mid rm = 0 \text{ for all } m \in M\}$  is
the trivial ideal~$\{0\}$. We call $M$ \emph{almost trivial} or \emph{non-faithful} if 
$M$ is not faithful, i.e.~there exists $r \in R, r \not= 0$ such that
$rm = 0$ holds for all $m \in M$. Recall for an $R$-module $M$ that $\tors(M)$ is the submodule
consisting of elements $m \in M$ for which there exists some $r \in  R, r \not= 0$ (depending on
$m$) with $rm = 0$.
We call $M$ \emph{almost finitely generated} if $\tors(M)$ is almost trivial and
$M/\tors(M)$ is a finitely generated $R$-module. 
We will often abbreviate $M/\tors(M)$ by $M/\tors$.

Recall that a \emph{Serre subcategory} of an abelian category is an
abelian subcategory closed under subobjects, quotients, and extensions
(see for example~\cite[Exercise~10.3.2 on page~384]{Weibel(1994)}).

\begin{lemma}
\label{lem: almost trivial and almost finitely generated modules}
\begin{enumerate}
\item \label{lem: almost trivial and almost finitely generated modules: almost trivial}
  The full subcategory of~$R\text{-}\MODULES$ whose objects are the almost
  trivial
  modules is a Serre subcategory.
\item \label{lem: almost trivial and almost finitely generated modules: almost f.g}
  If $R$ is Noetherian then
  the full subcategory of~$R\text{-}\MODULES$ whose objects are the almost
  finitely generated
  modules is a Serre subcategory.
\end{enumerate}
\end{lemma}

\begin{proof}
\ref{lem: almost trivial and almost finitely generated modules: almost trivial}
Obviously a quotient module and a submodule of an almost trivial module are 
almost trivial again.  Let $M_0 \to M_1 \to M_2$ be an exact sequence of $R$-modules
such $M_0$ and $M_2$ are almost trivial. Choose $r_i$ in $\Ann(M_i)$ with $r_i \not = 0$ for
$i = 0,2$. One easily checks that the product $r_0r_2$ lies in
$\Ann(M_1)$ and is different from zero. Hence $M_1$ is almost trivial. 
\\[2mm]
\ref{lem: almost trivial and almost finitely generated modules: almost f.g}
Consider an exact sequence of $R$-modules
$0 \to M_0 \to M_1 \to M_2 \to 0$. It suffices to prove that
both $M_0$ and $M_2$ are almost finitely generated if and only if $M_1$ is almost finitely
generated. Define $C$ to be the cokernel of the induced map
$\tors(M_1) \to \tors(M_2)$ and $K$ to be the kernel of the induced map
$M_1/\tors \to M_2/\tors$. Then we get exact sequences
$$0 \to \tors(M_0) \to \tors(M_1) \to \tors(M_2) \to C \to 0,$$
$$\qquad\qquad 0 \to M_0/\tors \to K \to C \to 0, \qquad\text{and}$$
$$0 \to K \to M_1/\tors \to M_2/\tors \to 0.$$

Suppose that $M_1$ is almost finitely generated, 
i.e.~$\tors(M_1)$ is almost trivial
and $M_1/\tors $ is finitely generated. 
As $R$ is Noetherian, $K$ is finitely generated. Therefore $C$ is finitely generated
and $C = \tors(C)$, hence $C$ is almost trivial. We conclude from 
assertion~\ref{lem: almost trivial and almost finitely generated modules: almost trivial}
that $\tors(M_0)$ and $\tors(M_2)$ are almost trivial. Since $R$ is Noetherian and
$M_1/\tors $ is finitely generated, both $M_0/\tors $ and $M_2/\tors $
are finitely generated. This shows that both $M_0$ and $M_2$ are almost finitely
generated.

Suppose that both $M_0$ and $M_2$ are almost finitely generated.
We conclude from  assertion~\ref{lem: almost trivial and almost finitely generated modules: almost trivial}
that $C$ and $\tors(M_1)$ are almost trivial. We obtain an exact sequence
\begin{multline*} \hom_R(K,M_0/\tors ) \longrightarrow \hom_R(M_0/\tors ,M_0/\tors )
 \longrightarrow \Ext^1_R(C,M_0/\tors ). 
\end{multline*}
Since $C$ is almost trivial, we can find $r \in R, r \not=0$ which annihilates $C$ and hence
also $\Ext^1_R(C,M_0/\tors )$. Therefore 
$r \cdot \id\colon M_0/\tors \to M_0/\tors $ is sent to the zero element in
$\Ext^1_R(C,M_0/\tors )$. Hence we can find a $R$-map $f \colon K \to M_0/\tors $
whose restriction to $M_0/\tors $ is $r \cdot \id\colon M_0/\tors \to
M_0/\tors $. Consider $x \in \ker(f)$. Since $r$ annihilates $C$, the element $r \cdot x$ lies in $M_0/\tors $
and belongs to the kernel of $r \cdot \id\colon M_0/\tors \to
M_0/\tors $. This implies $r^2 \cdot x = 0$. Since $K$ is a submodule of
$M_1/\tors $ and hence is torsionfree, we conclude $x = 0$. This shows
that $f \colon K \to M_0/\tors $ is an injection. Since $M_0/\tors $ 
is finitely generated and $R$ is Noetherian, $K$ is finitely generated. 
Since $M_2/\tors $  is finitely generated and $R$ is Noetherian, $M_1/\tors $ 
is finitely generated. Hence $M_1$ is almost finitely generated.
\end{proof}

\begin{definition} \label{def: almost isomorphism}
A map $f \colon M \to N$ of $R$-modules is called an \emph{almost
  isomorphism} if its kernel and cokernel are almost trivial
  $R$-modules.
\end{definition}

We conclude from Lemma~\ref{lem: almost trivial and almost finitely generated modules}~\ref{lem: almost trivial and almost finitely generated modules: almost trivial}
that for two composable $R$-maps $f \colon L \to M$ and $g \colon M
\to N$ all three maps $f$, $g$ and $g \circ f$ are almost
isomorphisms if two of them are. Moreover there is a Five-Lemma for
almost isomorphisms. For an  almost $R$-isomorphism $f$ the $F$-map
$F \otimes_R f$ is an $F$-isomorphism, where $F$ denotes the quotient
field of $R$. An $R$-module $M$ is almost trivial or almost finitely generated respectively
if and only if it is almost isomorphic to the zero module or respectively to a finitely
generated $R$-module.

\typeout{--------------------  Commuting holim and smash  --------------------}

\section{Commuting homotopy limits and smash products}
\label{sec: Commuting homotopy limits and smash products}

In this section we will prove our main result Theorem~\ref{the: commuting wedge and homotopy limit}
and the Addendum~\ref{addendum}. The proof will proceed by induction over the skeleta of $Y$ and hence we start with the 
following lemma.

\begin{lemma} \label{the: commuting wedge and homotopy limit: Y zero-dimensional}
Theorem~\ref{the: commuting wedge and homotopy limit} and Addendum~\ref{addendum} 
are true if $Y$ is a zero-dimensional $\cali$-$CW$-complex.
\end{lemma} 
\begin{proof}
For a large part we will treat the proofs of the theorem and of the addendum simultaneously.

By assumption we can find a family of objects $\{c_i \mid i \in I\}$ such that
$Y = \coprod_{i \in I} \mor_{\cali}(?,c_i)$. Notice that
there is a canonical isomorphism
\[
\map_{\cali}(\mor_{\cali}(?,c)_+,\bfF) \xrightarrow{\cong} \bfF(c)
\]
for any contravariant $\cali$-spectrum $\bfF$ and any object~$c$ in~$\cali$. 
Hence it suffices to show that the map
\[
\bft\colon   X_+ \sma_{\Or(G)} \prod_{i \in I} \bfE(c_i) \to 
\prod_{i \in I} \bfR\left(X_+ \sma_{\Or(G)} \bfE(c_i)\right)
\]
which is induced by the various projections $\pr_i\colon  \prod_{i \in I} \bfE(c_i) \to  \bfE(c_i)$ followed by the 
replacement maps 
$\bfr(X_+ \sma_{\Or(G)} \bfE(c_i)) \colon  X_+ \sma_{\Or(G)} \bfE(c_i) \to 
\bfR\left(X_+ \sma_{\Or(G)} \bfE(c_i)\right)$
 of~\eqref{bfr}, 
induces an (almost) isomorphism on homotopy groups $\pi_p$ for all $p \le n$. 

For all $p$ there is a commutative diagram
\begin{equation} \label{product-square}
\vcenter{\xymatrix{
\pi_p \left( X_+ \sma_{\Or (G) } \prod_{i \in I} \bfE (c_i) \right) \ar[r]^-{\pi_p ( \bft )} \ar[d] &
\pi_p \left(  \prod_{i \in I} \bfR \left( X_+ \sma_{\Or (G) } \bfE (c_i) \right) \right) \ar[d]^-{\cong} \\
\prod_{i \in I} \pi_p \left( X_+ \sma_{\Or (G) }  \bfE (c_i) \right) \ar[r]^-{\cong}  &
\prod_{i \in I} \pi_p \left( \bfR \left( X_+ \sma_{\Or (G) } \bfE (c_i) \right) \right) .
}}
\end{equation}
The vertical maps are induced by the obvious projections. The right hand \linebreak vertical map is an
isomorphism by Lemma~\ref{lem: commuting pi and products for Omega-spectra}. The lower horizontal map is the product of the isomorphisms 
induced by the fibrant replacement equivalences
$\bfr(X_+ \sma_{\Or(G)} \bfE(c_i))
\colon  X_+ \sma_{\Or(G)} \bfE(c_i) \to \bfR\left(X_+ \sma_{\Or(G)} \bfE(c_i)\right)$
and hence an isomorphism, compare~\eqref{bfr}.
We see that it suffices to show that the left hand vertical map
is an (almost) isomorphism 
for $p \le n$. This is what we will prove next.

For any covariant $\Or(G)$-spectrum $\bfF$ and any
 $G$-$CW$-complex $Z$ we get a homological spectral sequence
$(E_{*,*}^*,d_{*,*}^*)$ 
which converges to $\pi_{p+q}\left(Z_+ \sma_{\Or(G)} \bfF\right)$ and 
whose $E^2$-term is given by the Bredon homology
$$E_{p,q}^2 ~ = ~ H^{\bbZ \Or(G)}_p(Z;\pi_q(\bfF))$$
This follows from~\cite[Theorem~4.7 on page~234, Theorem~7.4~(3) on page~250]{Davis-Lueck(1998)}.
Recall that the \emph{Bredon homology} of a $G$-$CW$-complex $Z$ with
coefficients in a covariant $\bbZ\Or(G)$-module $M$ is defined by
\[
H_p^{\bbZ\Or(G)}(Z;M) = H_p\left(C_*^{\bbZ\Or(G)}(Z) \otimes_{\bbZ \Or(G)} M\right).
\]
This spectral sequence reduces to the classical Atiyah-Hirzebruch spectral sequence if $G$ is the trivial group.
A map of covariant $\Or(G)$-spectra $\bfF \to \bfF'$ induces a map of the associated spectral sequences.
Let $(E_{*,*}^*,d_{*,*}^*)$ be the spectral sequence associated to 
$\bfF = \prod_{i \in I} \bfE(c_i)$ and $(E[i]_{*,*}^*,d[i]_{*,*}^*)$ be the spectral sequence
associated to $\bfF = \bfE(c_i)$.
Our assumptions imply that all these spectral sequences live in the first quadrant (possibly shifted down by N).
This implies that all convergence questions are trivial, i.e. $E^{\infty}_{p,q} = E^r_{p,q}$ for
$r > p+q-N+1$ and the filtration of $\pi_{p+q}(X_+ \sma_{\Or(G)} \bfE(c_i))$, whose quotients are 
given by the $E^{\infty}$-term, consists of at most $p+q-N+1$ stages. 
Notice that $\prod_{i \in I}$ preserves exact sequences. Hence we can take the product
$\prod_{ i \in I} (E[i]_{*,*}^*,d[i]_{*,*}^*)$ and obtain a spectral sequence
which converges to 
$\prod_{i \in I} \pi_{p+q}\left(X_+ \sma_{\Or(G)} \bfE(c_i)\right)$ and whose $E^2$-term is
$\prod_{i \in I} H^{\bbZ \Or(G)}_p(X;\pi_q(\bfF))$.
The usual spectral sequence comparison argument shows that the left hand vertical map in~(\ref{product-square})
which is induced by the 
various projections $\pr_i\colon  \prod_{i \in I} \bfE(c_i) \to  \bfE(c_i)$ 
is an (almost) isomorphism for each $p \le n$, if the maps 
induced by the projections $\pr_i\colon  \prod_{i \in I} \bfE(c_i) \to  \bfE(c_i)$ 
\[
H_p^{\bbZ \Or(G)}\left(X_+;\pi_q\left(\prod_{i \in I} \bfE(c_i)\right)\right) 
\to \prod_{i \in I} H_p^{\bbZ \Or(G)}\left(X_+;\pi_q(\bfE(c_i))\right)
\]
are (almost) isomorphisms for all $p,q \in \bbZ$ with $p + q \le n$.

The canonical map
\[
\pi_q \left( \prod_{i \in I} \bfE(c_i) \right) \xrightarrow{\cong}  \prod_{i \in I} \pi_q(\bfE(c_i))
\]
is bijective by Lemma~\ref{lem: commuting pi and products for Omega-spectra} since we assume that the $\bfE (c_i)$ are
$\Omega$-spectra.
The assumptions about $\bfE$ imply that Lemma~\ref{lemma-factorize-over-groupoids} is 
applicable and hence
the $\bbZ\Or(G)$-modules
$\pi_q(\bfE(c_i))$ and $\prod_{i \in I} \pi_q(\bfE(c_i))$ factorize through the canonical projection
$\pr\colon  \Or(G) \to \Sub(G)$.
Let $\calf_X = \{H \subset G \mid X^H \not= \emptyset\}$
be the family of isotropy groups of $X$.
Let $\Sub(G;\calf_X)$ be the full subcategory 
$\Sub(G)$ whose objects are in~$\calf_X$.
Let $C^{\bbZ\Sub(G;\calf_X)}_*(X)$ be the contravariant
$\bbZ\Sub(G;\calf_X)$-chain complex
which is the composition of the ``cellular chain complex''-functor with the
contravariant functor from $\Sub(G;\calf_X)$ to the category of
$CW$-complexes which sends an object~$H$
to the $CW$-complex $Z_GH\backslash X^H$. 
Given a covariant $\bbZ\Sub(G;\calf_X)$-module~$N$, define 
$$H_p^{\bbZ \Sub(G;\calf_X)}(X;N) = 
H_p\left(C_*^{\bbZ\Sub(G;\calf_X)}(X) \otimes_{\bbZ\Sub(G;\calf_X)} N\right).$$
If $M$ is a covariant $\bbZ\Sub(G)$-module, then we obtain a natural isomorphism
\[
H_p^{\bbZ \Or(G)}(X;\pr^* M) \xrightarrow{\cong}
H_p^{\bbZ \Sub(G;\calf_X)}(X;M|_{\Sub(G;\calf_X)}),
\]
where $\pr^*M$ is the $\bbZ\Or(G)$-module $M \circ \pr$ and
$M|_{\Sub(G;\calf_X)}$ is the $\bbZ\Sub(G;\calf_X)$-module given by
the restriction of $M$ to $\Sub(G;\calf_X)$.
This follows from the adjunction of induction and
restriction~\cite[9.21 on page~169]{Lueck(1989)},
since the induction of $C_*^{\bbZ\Sub(G;\calf_X)}(X)$ with respect to the
inclusion $\Sub(G;\calf_X) \to \Sub(G)$
and the induction of $C_*^{\bbZ\Or(G)}(X)$ with respect to the projection
$\pr \colon \Or(G) \to \Sub(G)$ agree.
In the sequel we abbreviate $\Sub = \Sub(G;\calf_X)$. Hence it remains to
prove that for any family $\{M_i\mid i \in I\}$ of 
covariant $\bbZ \Sub$-modules the canonical map induced
by the various projections $\prod_{i \in I} M_i \to M_i$ 
$$H_p^{\bbZ \Sub}\left(X;\prod_{i \in I} M_i\right) ~ \to ~ \prod_{i \in I}  
H_p^{\bbZ \Sub}\left(X;M_i\right)$$
is an (almost) isomorphism  for all $p \in \bbZ$ with $p \le n - N$. 
We will prove that for a contravariant projective $\bbZ\Sub$-chain
complex $C_*$ (such as $C_*^{\bbZ\Sub(G;\calf_X)}(X)$) and a family of covariant $\bbZ\Sub$-modules $M_i$ 
the canonical map 
\begin{equation} \label{abzweigung}
H_p\left(C_* \otimes_{\bbZ\Sub} \prod_{i \in I} M_i\right) ~ \to ~ \prod_{i \in I}  
H_p\left(C_* \otimes_{\bbZ\Sub} M_i\right)
\end{equation}
induced
by the various projections $\prod_{i \in I} M_i \to M_i$ 
is an (almost) $\bbZ$-isomorphism for all $p \le n - N$, provided that the $\bbZ$-module
$H_p(C_*(H))$ is (almost) finitely generated
 for any object $H$ in $\Sub$ and $p \le n - N$.

We will first treat the case of the addendum, where we assume that the homology modules
are almost finitely generated and want to conclude an almost isomorphism.

There are spectral sequences converging to 
\[
H_{p+q}\left(C_* \otimes_{\bbZ\Sub} \prod_{i \in I} M_i\right) 
\quad \mbox{ respectively } \quad 
H_{p+q}\left(C_* \otimes_{\bbZ\Sub} M_i\right)
\]
whose $E^2_{p,q}$-term are given by 
\[
\Tor_p^{\bbZ \Sub} \left( H_q(C_*),\prod_{i \in I} M_i \right) 
\quad  \mbox{ respectively } \quad 
\Tor_p^{\bbZ \Sub}(H_q(C_*),M_i).
\]
Taking the product over $i \in I$ of the second 
spectral sequences yields a spectral sequence which converges to 
$\prod_{i \in I}  H_{p+q}\left(C_* \otimes_{\bbZ\Sub} M_i\right)$ and
whose $E^2_{p,q}$-term is 
$\prod_{i \in I}  \Tor_p^{\bbZ\Sub}(H_q(C_*),M_i)$.
By  the standard spectral sequence comparison
argument it suffices to prove that the canonical map 
\begin{eqnarray}
\Tor_p^{\bbZ\Sub}\left(H_q(C_*),\prod_{i \in I} M_i\right) \longrightarrow 
\prod_{i \in I}  \Tor_p^{\bbZ\Sub}(H_q(C_*),M_i)
\label{tor and products}
\end{eqnarray}
is an almost isomorphism for each $p,q \in \bbZ$ with  $p +q \le n -N$. 
The following diagram commutes
\comsquare{\Tor_p^{\bbZ\Sub}\left(H_q(C_*),\prod_{i \in I} M_i\right)}
{} 
{\prod_{i \in I}  \Tor_p^{\bbZ\Sub}(H_q(C_*),M_i)}
{}{}
{\Tor_p^{\bbZ \Sub}\left(H_q(C_*)/\tors ,\prod_{i \in I} M_i\right)}
{} 
{\prod_{i \in I}  \Tor_p^{\bbZ\Sub}(H_q(C_*)/\tors,M_i) .}
As $\Sub$ has only finitely many isomorphism classes of objects and 
$\tors(H_q(C_*(H)))$ is almost trivial for each object~$H$, there is an integer $k \not=0$
which annihilates $\tors(H_q(C_*(H)))$ for each object~$H$. Hence 
$\Tor_p^{\bbZ\Sub}(\tors(H_q(C_*)),N)$ is annihilated by $k$ for each covariant
$\bbZ\Sub$-module $N$. By the long exact $\Tor$-sequence associated to
the exact sequence of contravariant $\bbZ\Sub$-modules
\[
0 \to \tors(H_q(C_*)) \to H_q(C_*) \to H_q(C_*)/\tors(H_q(C_*)) \to 0
\]
we conclude
that both the kernel and the cokernel of the canonical map 
$$\Tor_p^{\bbZ \Sub}(H_q(C_*),N) \to \Tor_p^{\bbZ \Sub}(H_q(C_*)/\tors ,N)
$$
are annihilated by $k$ for any covariant $\bbZ\Sub$-module $N$. This implies that
in the commutative square above, the vertical arrows are almost isomorphisms.
Hence it suffices to prove that the lower
horizontal arrow is an almost $\bbZ$-isomorphism.
By assumption the $\bbZ$-module $H_q(C_*(H))/\tors $ is finitely generated for
all $q \le n - N$. 

Since we assume that $\Sub$ has only finitely many isomorphism classes of
objects and $\aut_{\Sub}(H) = W_GH = N_GH/H \cdot Z_GH$ is a finite group for
any object $H \in \Sub$, the category
$\Sub$ is finite in the sense of~\cite[Definition~16.1 on page~325]{Lueck(1989)}.
(Notice that $\Or(G;\calf_X)$ does not have this property and therefore
it is crucial to pass to $\Sub$.) 
Hence the category of $\bbZ\Sub$-modules is noetherian, in the sense that
each submodule of a finitely generated $\bbZ\Sub$-module
is itself finitely generated.  Moreover a covariant
$\bbZ\Sub$-module $L$ is finitely generated if and only if $L(H)$ is a finitely
generated $\bbZ$-module for each object $H \in \Sub$
(see~\cite[Lemma~16.10 on page~327]{Lueck(1989)}). By assumption 
$H_q(C_*(H))/\tors$ is a finitely
generated $\bbZ$-module for each object $H$ in $\Sub$. 
This implies that there is a finitely generated free (not necessarily finite-dimensional)
$\bbZ\Sub$-resolution $F_*$ for $H_q(C_*)/\tors $ for $q \le n - N$
(see~\cite[Lemma~17.1 on page~339]{Lueck(1989)}). Thus we get 
for $q \le n - N$ identifications
\[
\Tor_p^{\bbZ\Sub}(H_q(C_*)/\tors ,\prod_{i \in I} M_i) ~ = ~ 
H_p(F_* \otimes_{\bbZ\Sub} \prod_{i \in I} M_i)
\]
and
\begin{eqnarray*}
\prod_{i \in I} \Tor_p^{\bbZ\Sub}(H_q(C_*)/\tors ,M_i) & = & 
\prod_{i \in I} H_p(F_* \otimes_{\bbZ\Sub}M_i)  \\
&  = &  
H_p\left(\prod_{i \in I} F_* \otimes_{\bbZ\Sub} M_i\right),
\end{eqnarray*}
under which the map 
$$\Tor_p^{\bbZ \Sub}(H_q(C_*),\prod_{i \in I} M_i) \to 
\prod_{i \in I}  \Tor_p^{\bbZ \Sub}(H_q(C_*),M_i)$$
becomes the map induced on homology by the canonical $\bbZ$-chain map
\begin{equation} \label{final-step}
F_* \otimes_{\bbZ\Sub} \prod_{i \in I} M_i \to 
\prod_{i \in I} F_* \otimes_{\bbZ\Sub}  M_i.
\end{equation}
But this chain map is a chain isomorphism since each chain module
$F_i$ is a finitely generated free $\bbZ\Sub$-module. Here one uses the fact that the processes
of taking sums over a finite index set and products over
an arbitrary index set commute.

The case where we assume the homology groups to be finitely generated $\bbZ$-modules is similar but easier.
Using the fact that the category of $\bbZ\Sub$-modules is noetherian
one can replace the complex $C_{\ast}$ which appears in~(\ref{abzweigung})
by a homotopy equivalent complex~$F_{\ast}$ of projective $\bbZ\Sub$-modules which is of finite type. For such a complex
a map like the one in~(\ref{final-step}) is an isomorphism and 
hence so is the map~(\ref{abzweigung}).
\end{proof}

\begin{remark} \label{rem: Q-finitely generated too weak}
Before we finish the proof of Theorem~\ref{the: commuting wedge and homotopy limit} and its addendum 
we want to discuss where 
the proof above breaks down if we ask for a rational isomorphisms and 
weaken assumption (D) to the requirement 
that $H_n(Z_GH\backslash X^H;\bbQ)$ is a finitely
generated $\bbQ$-module for all $n$.

We would have to show
that the map appearing in~\eqref{tor and products} is a rational isomorphism.
The following example shows that this is not possible. 
We claim that for a fixed prime $p$ the canonical map
$$\Tor_1^{\bbZ}\left(\bigoplus_{m \ge 2} \bbZ/p^m,\prod_{n \ge 2} \bbZ/p^n\right)
~ \to ~
\prod_{n \ge 2}  \Tor_1^{\bbZ}\left(\bigoplus_{m \ge 2} \bbZ/p^m,\bbZ/p^n\right)
$$
cannot be a rational isomorphism although
$( \bigoplus_{m \ge 2} \bbZ/p^m ) \otimes_{\bbZ} \bbQ = 0$. 
Namely, the source is
$\bigoplus_{m \ge 2} \prod_{n \ge 2} \bbZ/p^{\min(m,n)}$, which vanishes after applying
$- \otimes_{\bbZ} \bbQ$. The target is
$\prod_{n \ge 2} \bigoplus_{m \ge 2}  \bbZ/p^{\min(m,n)}$, which contains
$\prod_{n \ge 2} \bbZ/p^n$ as a submodule and hence does not vanish after applying
$- \otimes_{\bbZ} \bbQ$. 
\end{remark}

Now we finish the proof of Theorem~\ref{the: commuting wedge and homotopy limit} and Addendum~\ref{addendum}.

\begin{proof}[Proof  of Theorem~\ref{the: commuting wedge and homotopy limit} and Addendum~\ref{addendum}.]
We use induction over the dimension $d$ of $Y$.  If $d = 0$, the claim has already been proved in
Lemma~\ref{the: commuting wedge and homotopy limit: Y zero-dimensional}.
The induction step from $d-1$ to $d \ge 1$ is done as follows. We can write $Y$ as a pushout
\comsquare{Y_0}{i_1}{Y_1}{i_2}{j_1}{Y_2}{i_2}{Y}
where $Y_1$ is a $(d-1)$-dimensional $\cali$-$CW$-complex, 
$Y_0 = \coprod_{i \in I} \mor_{\cali}(?,c_i) \times S^{d-1}$,
$Y_2 = \coprod_{i \in I} \mor_{\cali}(?,c_i) \times D^d$ and $i_2$ is the inclusion. 
In particular $Y_1$ is a $(d-1)$-dimensional $CW$-complex and 
$Y_2$ is $\cali$-homotopy equivalent to the $0$-dimensional $CW$-complex
$\coprod_{i \in I} \mor_{\cali}(?,c_i)$. Hence the induction hypothesis applies to
$Y_0$, $Y_1$ and $Y_2$.

The canonical map from the homotopy pushout of $Y_2  \leftarrow Y_0 \rightarrow Y_1$ to $Y$ 
is a map of $\cali$-$CW$-complexes whose evaluation at each object $c$ is a homotopy equivalence
since $i_2(c)$ is a cofibration. Hence it is a $\cali$-homotopy equivalence
(see~\cite[Corollary~3.5 on page~222]{Davis-Lueck(1998)}).
Since 
for any covariant $\cali$-space $Z$ 
the functor $\map_{\cali}(-,Z)$ sends 
a homotopy cocartesian square of spaces to a homotopy cartesian one,
the following diagram is levelwise  homotopy cartesian for each object $G/H \in \Or(G)$
\comsquare{\map_{\cali}(Y_+,\bfE(G/H))}{}{\map_{\cali}((Y_1)_+,\bfE(G/H))}
{}{}
{\map_{\cali}((Y_2)_+,\bfE(G/H))}{}{\map_{\cali}((Y_0)_+,\bfE(G/H)).}
We conclude with Lemma~\ref{lem: hopushout and hopullback for spectra} that the square is homotopy cocartesian.
By the same line of argument the diagram 
\comsquare{\map_{\cali}(Y_+,\bfR(X_+ \sma_{\Or(G)}\bfE(G/?)))}{}
{\map_{\cali}((Y_1)_+,\bfR(X_+  \sma_{\Or(G)} \bfE(G/?)))}
{}{}
{\map_{\cali}((Y_2)_+,\bfR(X_+  \sma_{\Or(G)} \bfE(G/?)))}{}
{\map_{\cali}((Y_0)_+,\bfR(X_+  \sma_{\Or(G)} \bfE(G/?))) }
is homotopy cocartesian.
Since the functor $X_+  \sma_{\Or(G)} -$ respects homotopy cocartesian squares
because of the associativity of $\sma$,
we conclude from 
Lemma~\ref{lem: weak equivalences and compatibility with otimes and map}~\ref{lem: weak equivalences and compatibility with otimes and map: spectra and wedge}
that the following diagram is homotopy cartesian
\comsquare{X_+ \sma_{\Or(G)} \map_{\cali}(Y_+,\bfE)}{}
{X_+ \sma_{\Or(G)} \map_{\cali}((Y_1)_+,\bfE)}
{}{}
{X_+ \sma_{\Or(G)} \map_{\cali}((Y_2)_+,\bfE)}{}
{X_+ \sma_{\Or(G)} \map_{\cali}((Y_0)_+,\bfE) .}
The various transformations $\bft(X,Y_i,\bfE)$ and $\bft(X,Y,\bfE)$
yield a map from this homotopy cartesian diagram to the one before.
Hence they yield a transformation from the long exact homotopy sequences
associated in Lemma~\ref{lem: hopushout and hopullback for spectra}
to these homotopy pullbacks.
Now Theorem~\ref{the: commuting wedge and homotopy limit}
follows from the Five-Lemma for (almost) isomorphisms.
\end{proof}

\typeout{--------------------  The linear version  --------------------}

\section{The linear version}
\label{sec: The linear version}

In this section we state the linear version of our main
Theorem~\ref{the: commuting wedge and homotopy limit}.
This will be needed in~\cite{Lueck-Matthey-Reich(2003)}.
The proof in the linear setting is completely analogous to the one for
spaces and spectra presented in Section~\ref{sec: Commuting homotopy
limits and smash products}, and therefore we omit it.

\begin{theorem}\label{the: commuting otimes and hom}
Consider the following data:
\begin{itemize}
\item
A functor $E_* \colon \cali^{\op} \times \Or(G) \to \ch(\bbZ\text{-}\MODULES)$
which factorizes as
\[
E_* = F_* \circ (\id \times \calg^G),
\]
where $F_* \colon \cali^{\op} \times \GROUPOIDS \to \ch(\bbZ\text{-}\MODULES)$
has the property that for each object $c$ in~$\cali$
the functor $F_*(c, - )$ sends equivalences of groupoids to homotopy equivalences.
\item
A contravariant chain complex $D_*$ of projective $\bbZ\cali$-modules.
\item
A $G$-CW complex~$X$.
\end{itemize}
Suppose that there are numbers $d$, $n$ and~$N \in \bbZ$ with~$d \geq 0$
such that the following conditions are satisfied:
\begin{enumerate}
\item[(A)]
$D_*$ is $d$-dimensional, i.e.~$D_k(c)=0$ for all objects $c$ in~$\cali$
and for all $k\in\bbZ$ such that $k<0$ or~$k>d$.
\item[(B)]
The chain complexes $E_*(c , G/H)$ are uniformly homologically
$(N-1)$-bounded below,
i.e.~for all objects $c$ in~$\cali$ and all orbits~$G/H$
we have $H_q (E_*(c,G/H)) = 0 $ for~$q < N$.
\item[(C)]
The $G$-CW complex~$X$ has only finite isotropy groups
and only finitely many orbit types.
\item[(D)]
For each finite subgroup $H \subset G$ and every $p \leq n + d - N$
the homology group
\[
H_p ( Z_G H \backslash X^H ; \bbZ )
\]
is an (almost) finitely generated $\bbZ$-module.
\end{enumerate}
Then for each $p \leq n$ the map
\[
H_p\left( C_*^{\bbZ\Or(G)}(X) \otimes_{\bbZ\Or(G)} \hom_{\cali}(D_*,E_*) \right)
\xrightarrow{H_p(t_*)}
H_p\left( \hom_{\cali}\left(D_* , C_*^{\bbZ\Or(G)}(X) \otimes_{\bbZ\Or(G)} E_*\right) \right)
\]
is an (almost) isomorphism.
\end{theorem}

Here of course
\[
t_* \colon  C_*^{\bbZ\Or(G)}(X) \otimes_{\bbZ\Or(G)} \hom_{\cali}(D_*,E_*)
\to
\hom_{\cali}\left(D_*, C_*^{\bbZ\Or(G)}(X) \otimes_{\bbZ\Or(G)} E_* \right)
\]
is the chain map that sends $x \otimes \phi$ to the map
$D_* \to C_*^{\bbZ\Or(G)}(X) \otimes_{\bbZ\Or(G)} E_*$ given by
$y \mapsto x \otimes \phi(y)$.
There are also linear versions of Examples~\ref{exa: Adams result}
and~\ref{exa: necessary finiteness conditions on X}.

We mention that Theorem~\ref{the: commuting otimes and hom} above
remains true in a more general setting.  Namely we can replace the
ring of integers~$\bbZ$ by any commutative Noetherian ring~$R$, and
$C_*^{\bbZ\Or(G)}(X)$ by any contravariant $R\Or(G)$-chain complex
satisfying the following conditions: 
$C_*$ is a complex of projectives concentrated in non-negative degrees; 
the set $S(C_*) := \{(H) \mid S_{G/H}C_p \not= 0 \text{ for some } p \in \bbZ\}$ is finite;
$(H) \in S(C_*)$ implies that~$H$ is finite; and 
for any finite group $H \subset G$ and any $p \le n + d - N$ the homology group $H_p((\pr_*C_*)(H))$ is an (almost)
finitely generated $\bbZ$-module.

Here $(H)$ denotes the conjugacy class of a subgroup $H \subseteq G$
and $S_{G/H}$ is the splitting functor defined in~\cite[9.26 on page~170]{Lueck(1989)}; 
$pr_*C_*$ is the $R\Sub(G)$-chain complex obtained from $C_*$ by induction
(see~\cite[9.15 on page~166]{Lueck(1989)})) with the projection
$\pr\colon \Or(G) \to \Sub(G)$ defined in~\eqref{pr colon Or(G) to
  Sub(G)}.  

The situation above is a special case because 
$C_*^{\bbZ\Or(G)}(X)$ is a free (in the sense
of~\cite[9.16 on page~167 and page~356]{Lueck(1989)}) and hence a
projective $\bbZ\Or(G)$-chain complex, the set
$S(C_*^{\bbZ\Or(G)}(X))$ is just the set of conjugacy classes of
isotropy groups of~$X$, and
$H_p(\pr_*C_*^{\bbZ\Or(G)}(X)(H))\cong H_p(Z_GH\backslash X^H;\bbZ)$.

\typeout{--------------------  Examples of categories  --------------------}

\section{Examples of categories with finite-dimensional classifying spaces}
\label{sec: Examples of categories with finite-dimensional classifying space}

In this section we show two examples of categories that have a
finite-dimensional contravariant classifying space, in the sense of
Definition~\ref{def: Ecalc}.  The second example is relevant to the construction
of topological cyclic homology and will play a role in the
applications discussed at the end of the introduction.  In both
examples we will denote by~$\bbN$ the set of non-negative
integers~$\{0,1,2,\ldots\}$.

\begin{example}\label{ex:Ecaln}
Let $\caln$ be the category associated to the partially ordered
set~$\bbN$, i.e.~the category whose objects are the non-negative
integers~$\bbN$ and for which $\mor_{\caln}(i,j)$ is empty if $i > j$,
or consists of precisely one element if $i \le j$.  There is a
$1$-dimensional model for the contravariant classifying space $E\caln$.
In fact, define $E\caln(i)$ to be the interval
$[i,\infty)$.  The zero-skeleton is given by the intersection of
$[i,\infty)$ with $\bbZ$; the 1-cells are then obvious.
(Compare also~\cite[Example~3.9 on page~224]{Davis-Lueck(1998)}.)
\end{example}

\begin{example}
Let $\calrf$ be the following category. The set of objects is 
$\bbN=\{0,1,2, \ldots\}$. The set of morphisms from $m$ to $n$ is
$\{(i,j) \in \bbN \times \bbN \mid i+j =n-m\}$.  In particular there
are no morphisms from $m$ to $n$ if $m > n$.  The identity morphism of
$n$ is given by $(0,0)$.  The composition of $(i_0,j_0)\colon n_0 \to
n_1$ and $(i_1,j_1)\colon n_1 \to n_2$ is $(i_0 + i_1, j_0+j_1)\colon
n_0 \to n_2$.

Notice that a contravariant functor~$\bfE$ from $\calrf$ to the
category of spectra is the same as a sequence of spectra $\bfE_0$, $\bfE_1$,
$\bfE_2$, $\ldots$ together with a system of maps
\[
\xymatrix{
\ldots \ar@<.7ex>[r]^{\bfr_2} \ar@<-.7ex>[r]_{\bff_2} &
\bfE_2 \ar@<.7ex>[r]^{\bfr_1} \ar@<-.7ex>[r]_{\bff_1} &
\bfE_1 \ar@<.7ex>[r]^{\bfr_0} \ar@<-.7ex>[r]_{\bff_0} &
\bfE_0
}
\]
such that $\bfr_{i} \circ \bff_{i+1} = \bff_{i} \circ \bfr_{i+1}$ holds
for all~$i \in \bbN$.
Namely, put $\bfE_i = \bfE(i)$, $\bfr_i = \bfE((0,1)\colon  i \to (i+1))$
and $\bff_i = \bfE((1,0)\colon  i \to (i+1))$.
Topological cyclic homology is defined as the homotopy
limit of such a system of spectra, where the~$\bfE_n$'s are given by
the fixed points of cyclic $p$-groups acting on topological Hochschild
homology, and the structure maps are the so-called Restriction and
Frobenius maps (see~\cite{Boekstedt-Hsiang-Madsen(1993)}).

\begin{proposition} \label{the: Ecalc is two-dimensional}
There is a two-dimensional model for the contravariant classifying
space of $\calrf$.
\end{proposition}

\begin{proof}
Define a contravariant $\calrf$-space $E$ by sending an object $n$ to
the space $\{(x,y)\in\bbR\times\bbR \mid x,y \ge 0\}$ and a morphism
$(i,j)\colon  m \to n$ to the map $(x,y) \mapsto (x+i,y+j)$.
Obviously each space $E(n)$ is contractible.  It remains to 
construct the $\calrf$-$CW$-structure. 

Define $E_0(n) = \bbN \times \bbN$, 
$E_1(n) = \{(x,y) \in \bbR\times\bbR \mid x,y \ge 0, \text{$x$ or~$y$}\in
\bbN\}$ and $E_2(n) = E$. This yields a filtration 
$E_0 \subset E_1 \subset E_2 = E$ of contravariant $\calrf$-spaces.
We claim that this is an $\calrf$-$CW$-structure. We have to construct
the relevant pushouts of contravariant $\calrf$-spaces.
Define for each $n \in \bbN$ embeddings of contravariant $\calrf$-spaces
\begin{eqnarray*}
Q_0(n) \colon   \mor_{\calrf}(?,n) \to E
& & (i,j) \mapsto (i,j);
\\
Q_1^h(n) \colon   \mor_{\calrf}(?,n) \times [0,1] \to E
& & ((i,j),s)  \mapsto (i+s,j);
\\
Q_1^v(n) \colon   \mor_{\calrf}(?,n) \times [0,1] \to E
& & ((i,j),s)  \mapsto (i,j+s);
\\
Q_2(n) \colon   \mor_{\calrf}(?,n) \times [0,1]
\times[0,1] \to E
& & ((i,j),s,t )  \mapsto (i+s,j+t).
\end{eqnarray*}
We claim that $\coprod_{n \in \bbN} \im(Q_0) = E_0$ and 
that we get pushouts of contravariant $\calrf$-$CW$-complexes
\[
\xymatrix@C=8em{
\coprod_{n\in\bbN} \Bigl( \mor_{\calrf}(?,n) \times \{0,1\}
                  \coprod \mor_{\calrf}(?,n) \times \{0,1\} \Bigr)
\ar[r]^-{\coprod_{n\in\bbN}\left(Q_1^h| \coprod Q_1^v|\right)} \ar[d]
& E_0 \ar[d] \\
\coprod_{n\in\bbN} \Bigl( \mor_{\calrf}(?,n) \times  [0,1]
                  \coprod \mor_{\calrf}(?,n) \times  [0,1]  \Bigr)
\ar[r]^-{\coprod_{n\in\bbN}\left(Q_1^h  \coprod Q_1^v \right)}
& E_1
}
\]
and
\[
\xymatrix@C=5em{
\coprod_{n\in\bbN} \mor_{\calrf}(?,n) \times \partial\bigl([0,1]\times[0,1]\bigr)
\ar[r]^-{\coprod_{n\in\bbN} Q_2|} \ar[d]
& E_1 \ar[d] \\
\coprod_{n\in\bbN} \mor_{\calrf}(?,n) \times               [0,1]\times[0,1]
\ar[r]^-{\coprod_{n\in\bbN} Q_2}
& E
}
\]
where $Q_1^h|$, $Q_1^v|$ and $Q_2|$ denote the obvious restrictions.
This finishes the proof of Proposition~\ref{the: Ecalc is two-dimensional}.
\end{proof}
\end{example}

\typeout{--------------------  The homology of EG(F)  --------------------}

\section{The homology of classifying spaces for families}
\label{sec: The homology of classifying spaces for families}

This section is devoted to the proof of 
Proposition~\ref{cor: assumptions of main theorem}.
We first recall some definitions.

A \emph{$G$-$CW$-complex} is a $CW$-complex with $G$-action
such that for an open cell $e$ and $g \in G$ with $g\cdot e \cap e \not= \emptyset$
multiplication with $g$ induces the identity on $e$.
For information about $G$-$CW$-complexes we refer for instance
to~\cite[Sections~1 and~2]{Lueck(1989)}
and~\cite[Chapter~II, Section~1]{Dieck(1987)}. A $G$-CW complex is proper if all isotropy groups are finite.
Recall that a \emph{family} $\calf$ of subgroups of $G$ is a set of subgroups closed under
conjugation and taking subgroups.
A $G$-$CW$-complex $\EGF{G}{\calf}$ 
all whose isotropy groups belong to 
$\calf$ and whose $H$-fixed point sets are contractible whenever $H \in
\calf$ is called a \emph{classifying space for the family} $\calf$.
This notion is due to tom Dieck (see~\cite[Chapter~I, Section~6]{Dieck(1987)}).
Any $G$-$CW$-complex whose isotropy groups belong to
$\calf$
admits up to $G$-homotopy precisely one $G$-map to $\EGF{G}{\calf}$. 
In particular any two models for $\EGF{G}{\calf}$ are $G$-homotopy equivalent. A
functorial bar-type construction for $\EGF{G}{\calf}$ is given in~\cite[Lemma~7.6]{Davis-Lueck(1998)}.

\begin{lemma} \label{lem: map on homology induced by EG times X to G  backslash X}
Let $X$ be a proper $G$-$CW$-complex. Suppose that there exists for
each $p \ge 0$ a   positive integer $d(p)$ such that for any isotropy group
$H$ of $X$ multiplication with $d(p)$ annihilates $\widetilde{H}_p(BH;\bbZ)$.

Then the canonical projection 
$\pr_X\colon EG \times_G X \to G\backslash X$  induces for all~$p \ge 0$
a $\bbZ$-almost isomorphism
$$H_p(\pr_X;\bbZ) \colon H_p( EG \times_G X ;\bbZ) \to H_p(G\backslash X;\bbZ).$$
\end{lemma}
\begin{proof}
Let $X_n$ be the $n$-skeleton of $X$. Then both maps $EG \times_G X_n
\to EG \times_G X$ and $G\backslash X_n \to G\backslash X$ are
$n$-connected and induce isomorphisms on $H_p(-;\bbZ)$ for $p < n$. Hence
it suffices to prove the claim for each $n$-skeleton $X_n$. This will
be done by induction over $n$. The induction beginning $n = -1$ is
trivial since $X_{-1} = \emptyset$. The induction step from $n$ to
$n+1$ is done as follows.

We can write $X_{n+1}$ as a $G$-pushout
\comsquare{\coprod_{i \in I} G/H_i \times S^n}{\coprod_{i \in I} q_i}{X_n}
 {}{}
{\coprod_{i \in I} G/H_i \times D^{n+1}}{\coprod_{i \in I}
  Q_i}{X_{n+1}.}
Applying $EG \times_G - $ and $G \backslash - $ one obtains two 
pushout squares whose left vertical arrows are again cofibrations.
The projections $\pr_Z$ where $Z$ runs through the four corners of the $G$-pushout square
yield a map between the two new squares. 
By a Five-Lemma argument, the K\"unneth
formula and the fact that $H_p(-;\bbZ)$ satisfies the disjoint union axiom for arbitrary
index sets, the induction step follows from the assumption that the map
$H_p(\pr_{G/H};\bbZ) \colon H_p(EG \times_G G/H;\bbZ) \to H_p(\{pt\};\bbZ)$ 
is surjective and its kernel is annihilated by $d(p)$ for each isotropy group $H$ of $X$.
\end{proof}

\begin{lemma} \label{lem: homological condition for the classifying space}
Let $\calf$ be a family of finite subgroups of $G$.
Fix an element $H \in \calf$. Suppose that there exists for
each $p \ge 0$ a positive integer $d(p)$ such that for any $K \in \calf$
with $H \subset K$ the annihilator ideal of $\widetilde{H}_p(B(Z_GH \cap K);\bbZ)$ contains $d(p)$.  Then
$H_p(Z_GH\backslash\EGF{G}{\calf}^H;\bbZ)$ is an almost finitely generated $\bbZ$-module for all $p \ge 1$.
\end{lemma}

\begin{proof}
Fix $H \in \calf$.  The projection
$EZ_GH \times_{Z_GH} \EGF{G}{\calf}^H  \to BZ_GH$ is a homotopy equivalence
since $\EGF{G}{\calf}^H$ is contractible. 
Now apply Lemma~\ref{lem: map on homology induced by EG times X to G  backslash X} to
the proper $Z_GH$-$CW$-complex $\EGF{G}{\calf}^H$ and 
the projection $EZ_GH \times_{Z_GH} \EGF{G}{\calf}^H \to Z_GH\backslash \EGF{G}{\calf}^H $.
\end{proof}

We are now prepared to prove Proposition~\ref{cor: assumptions of main theorem}.

\begin{proof}[Proof of Proposition~\ref{cor: assumptions of main theorem}]
Obviously $\EGF{G}{\calf}$  has finite isotropy groups and only finitely many orbit types.
There is a number $d$ such that the order of any element $H \in \calf$ divides~$d$.
Hence $\widetilde{H}_p(BK;\bbZ)$ is annihilated by~$d$ for each $K \in \calf$
and~$p$~\cite[Section~III.10]{Brown(1982)}.
Now apply Lemma~\ref{lem: homological condition for the classifying space}. 
\end{proof}

\typeout{--------------------  References  --------------------}



\end{document}